\documentclass[twocolumn]{autart}    

\usepackage{amsmath,amssymb,amsfonts,bm}
\usepackage{graphicx}
\usepackage{algorithm}
\usepackage{algorithmic}
\usepackage{cite}
\usepackage{color}
\let\theoremstyle\relax

\usepackage{amsthm}
\theoremstyle{plain}
\newtheorem{theorem}{Theorem}[section]
\newtheorem{lemma}[theorem]{Lemma}
\newtheorem{corollary}[theorem]{Corollary}
\newtheorem{proposition}[theorem]{Proposition}
\theoremstyle{definition}
\newtheorem{definition}[theorem]{Definition}
\theoremstyle{remark}
\newtheorem*{remark}{Remark}
\usepackage{wrapfig}

\begin{document}

\begin{frontmatter}

\title{A Non-Gaussian Bayesian Filter Using Power and Generalized Logarithmic Moments} 

\author[Guangyu]{Guangyu Wu}\ead{chinarustin@sjtu.edu.cn},    
\author[Anders]{Anders Lindquist}\ead{alq@kth.se}              

\address[Guangyu]{Department of Automation, Shanghai Jiao Tong University, Shanghai, China}  
\address[Anders]{School of Artificial Intelligence, Anhui University, Hefei, China and the Smart Sensor Fusion Research Laboratory, Shanghai Jiao Tong University, Shanghai, China}


\begin{abstract}                          
In this paper, we  propose a consistent non-Gaussian Bayesian filter for which the system state is a continuous function. The distributions of the true system states, and those of the system and observation noises, are only assumed Lebesgue integrable with no prior constraints on what function classes they fall within. This type of filter has significant merits in both theory and practice, in that it is able to ameliorate the curse of dimensionality present in the particle filter, a popular non-Gaussian Bayesian filter for which the system state is parameterized by discrete particles and  corresponding weights. We first propose a new type of statistics, called the generalized logarithmic moments. Together with the power moments, they are used to form a density surrogate, parameterized as an analytic function, to approximate the true system state. The map from the parameters of the proposed density surrogate to both the power moments and the generalized logarithmic moments is proved to be a diffeomorphism, establishing the fact that there exists a unique density surrogate which satisfies both moment conditions. This diffeomorphism also allows us to use gradient methods to treat the convex optimization problem in determining the parameters. Last but not least, simulation results reveal the advantage of using both sets of moments for estimating mixtures of complicated types of functions. A robot localization simulation is also given, as an engineering application to validate the proposed filtering scheme.
\end{abstract}

\end{frontmatter}

\section{Introduction}
\label{Intro}
In this paper, we consider the non-Gaussian Bayesian filtering problem for first-order systems following our previous work \cite{wu2023non}. The Bayesian filter offers a cohesive and recursive solution to the general stochastic filtering challenges. Ho and Lee's work \cite{ho1964bayesian} represents one of the initial attempts of iterative Bayesian estimation, delineating the principles and procedures of Bayesian filtering. Sprangins \cite{spragins1965note} delved into the iterative application of Bayes' rule for sequential parameter estimation. Lin and Yau \cite{lin1967bayesian}, as well as Chien and Fu \cite{chien1967bayesian}, explored the Bayesian approach for optimizing adaptive systems. The Bayesian filter consists of an iterative measurement-time update process, sometimes referred to by different terms. During the time update step, the system equation calculates a one-step ahead prediction of the state. In the measurement update step, the observation equation computes the correction to the state estimate according to the current observation.

The Bayesian filter usually does not provide an analytic state estimator, except for the cases where the system and observation equations are linear, and the distributions of the system state, the system noise and the observation noise are Gaussian. In the well-known Kalman filter (and its extended forms such as the extended Kalman filter and the unscented Kalman filter), to determine the integral in the time update step is essentially a parametric estimation problem, which is done by estimating the first and second order moments \cite{kalman1960new, kalman1961new, anderson1971kalman, kai2009robust, chen1993approximate, chen2000mixture}.

However, when the distributions are non-Gaussian, the problem is much more complicated. Given that the distributions of the system states and the noises are not conjugate pairs, the time update step does not provide an analytic probability density function of the system state, due to the intractable convolution \cite{chen2003bayesian}. The density function can then only be obtained by approximation. In the real-world applications, the distributions are usually non-Gaussian, which makes this approximation problem a significant one both in theory and in practice. Being an open problem, density approximation of the intractable integral has been a core problem of Bayesian filtering for decades and is still a hot topic. Numerous numerical methods have been proposed to obtain an analytic solution to the convolution in the time update step. To name a few, there are Gaussian/Laplace approximation \cite{mackay1998choice}, iterative quadrature \cite{freitas1999bayesian, kushner2000nonlinear, wang1978optimal}, Gaussian sum approximation \cite{sorenson1971recursive, alspach1972nonlinear}. In these methods, the intractable integral is approximated by a single Gaussian density function or a weighted sum of Gaussians. By doing this, the convolution operation in the time update step is tractable again, which makes it feasible for us to obtain the integral in an analytic form of function. A non-Gaussian Bayesian filter based on state-space calculus is proposed in \cite{hanzon2001state}. Instead of transforming the density function into a Gaussian or a mixture of Gaussians, this method considers the rational probability density functions, which are transfer functions of finite-dimensional linear systems by the realization theory. It provides explicit state space
descriptions for products and convolutions of rational densities, which then provides an analytic density function of the system state in a rational form. These parametric methods assume that the prior density belongs to a specified function class. This causes the flexibility of these methods to be limited. If the density function does not fall exactly within the assumed class, the estimation result may be severely biased. 

Since the problem we treat does not restrict the non-Gaussian density to fall within specific classes of functions, estimating the intractable prior density in the time update step is indeed an infinite-dimensional problem. The particle filter treats this estimation problem using discrete points without any assumption on the form of function of the prior density, which also turns the infinite-dimensional problem into a finite dimensional and tractable one \cite{carpenter1999improved, doucet2001sequential, doucet2001particle, andrieu2002particle, vermaak2002particle}. However characterizing the prior by discrete points requires massive particles to store the probability values of the states. The problem is even worse with the increase of dimensions, which is due to the curse of dimensionality \cite{chen2003bayesian}. Moreover, analyzing the errors of the particle filters is an extremely difficult task due to the indeterministic estimates caused by the sampling strategy \cite{Fox_2002, fox2003adapting}, and its performance suffers a lot from sample depletion \cite{gustafsson2010particle, doucet2001particle}. A non-Gaussian Bayesian filter, of which the system state distribution is a continuous function of a limited number of parameters, possesses notable advantages and is sought after by researchers in stochastic filtering.

In our previous papers \cite{wu2023non, wu2022multivariate}, a non-Gaussian Bayesian filter is implemented by approximating the intractable integral using the power moments. A non-classical density surrogate for the system state, in the form of a continuous function, is proposed, and the parameters of the proposed parametrization can also be determined by a convex optimization scheme with moment constraints, to which the solution is proved to exist and be unique. By the proposed filter, the power moments of the density estimates asymptotically converge to the true moments. Furthermore, according to Theorem 4.5.5 in \cite{chung2001course}, as the number of moment terms used approaches infinity, the approximated integral approaches the true analytic density in probability, which reveals the fact that the proposed filter is consistent. In the field of stochastic filtering, exact filters can only be computed for models that are completely discrete or for discrete-time linear Gaussian models, where the Kalman filter can be applied \cite{papanicolaou2014stochastic}. For system models that do not have exact filters, the best filters we can design for them are consistent filters. Therefore, the consistency of the filter proposed in \cite{wu2023non} is a clear advantage over other alternative filtering methods. We would also like to emphasize that an error upper bound in the sense of total variation distance exists for this filter. Such an error upper bound has not been proposed for other non-Gaussian Bayesian filters. With a relatively longer execution time for each filtering step, the proposed filter overcomes the disadvantages of the particle filters mentioned above. Moreover, the proposed filter can treat the filtering problem where the probability density of the true system state is an arbitrary Lebesgue-integrable function with first several orders of power moments being finite, without assuming the density to fall within specific function classes. 

In this paper, inspired by \cite{2002Identifiability}, we propose to use logarithmic-type moments together with power moments to improve the performance of non-Gaussian Bayesian filtering. The paper is organized as follows. In Section 2, we note that the conventional logarithmic moment does not work in this problem, and we propose a novel generalized logarithmic moment. An algorithm framework for Bayesian filtering using both power and generalized logarithmic moments is also proposed. Then we prove that by our proposed algorithm, the generalized logarithmic moments of the density estimates are asymptotically unbiased and approximately identical to the true ones. In Section 3, together with the fact that the power moments of the state estimate are approximately identical to the true ones given a large $n$, we propose to use both the power and generalized logarithmic moments to parameterize the density of the state. Then in Section 4, we prove that the parameters of the proposed density surrogate can be uniquely determined in terms of the power and generalized logarithmic moments up to order $2n$, by proving the corresponding map being diffeomorphic. Three density approximation examples are performed in Section 5, including the mixtures of Gaussian, generalized logistic and Laplacian densities. The simulation results with a comparison to the parametrization using only the power moments validate the advantage of using both sets of moments for parametrization of the prior density. Moreover, the proposed filtering scheme is applied to a robot localization task and the performance is compared to prevailing methods including the Kalman filter and the particle filter.

\section{Non-Gaussian Bayesian filtering and the generalized logarithmic moments}

Consider the stochastic system
\begin{equation}
\label{System}
    \begin{aligned}
    x_{t+1} &=f_{t} x_{t}+\eta_{t} \\
    y_{t} &=h_{t} x_{t}+\epsilon_{t}
\end{aligned}    
\end{equation}
$t=0,1,2, \ldots$, where the state $x_{t}$ is a random variable defined on $\mathbb{R}$, and $f_{t}, h_{t}$ are assumed to be known real numbers. The system noise $\eta_{t}$ and the observation noise $\epsilon_{t}$ are assumed to be random variables for which the distributions are Lebesgue integrable functions. Moreover, the noises are assumed to be independent from each other, and their distributions are denoted as $\rho_{\eta_{t}}$ and $\rho_{\epsilon_{t}}$. 

We use the Bayesian filtering framework in \cite{hanzon2001state}. Then the conditional probability density functions of the measurement and time updates are given by

\noindent \textbf{Measurement update}: For $t=0$,
\begin{equation}
\begin{aligned}
    \rho_{x_{0} \mid \mathcal{Y}_{0}}(x) = & \frac{\rho_{y_{0} \mid x_{0}}\left(y_{0}\right) \rho_{x_{0}}(x)}{\int_{\mathbb{R}}\rho_{y_{0} \mid x_{0}}\left(y_{0}\right) \rho_{x_{0}}(x)dx}\\
    = & \frac{\rho_{\epsilon_{0}}\left(y_{0}-h_{0} x\right) \rho_{x_{0}}(x)}{\int_{\mathbb{R}}\rho_{\epsilon_{0}}\left(y_{0}-h_{0} x\right) \rho_{x_{0}}(x)dx};
\label{Update1}
\end{aligned}
\end{equation}
for $t \geq 1$,
\begin{equation}
\begin{aligned}
\rho_{x_{t} \mid \mathcal{Y}_{t}}(x) = & \frac{\rho_{y_{t} \mid x_{t}}\left(y_{t}\right) \rho_{x_{t} \mid \mathcal{Y}_{t-1}}(x)}{\int_{\mathbb{R}}\rho_{y_{t} \mid x_{t}}\left(y_{t}\right) \rho_{x_{t} \mid \mathcal{Y}_{t-1}}(x)dx}\\
= & \frac{\rho_{\epsilon_{t}}\left(y_{t}-h_{t} x\right) \rho_{x_{t} \mid \mathcal{Y}_{t-1}}(x)}{\int_{\mathbb{R}}\rho_{\epsilon_{t}}\left(y_{t}-h_{t} x\right) \rho_{x_{t} \mid \mathcal{Y}_{t-1}}(x)dx}.
\label{Update2}
\end{aligned}
\end{equation}

\noindent \textbf{Time update}: For $t \geq 0$,
\begin{equation}
\begin{aligned}
    \rho_{x_{t+1} \mid \mathcal{Y}_{t}}(x) = & \left(\rho_{f_{t} x_{t} \mid \mathcal{Y}_{t}} * \rho_{\eta_{t}}\right)(x)\\
    = & \int_{\mathbb{R}} \rho_{x_{t} \mid \mathcal{Y}_{t}}\left(\frac{\varepsilon}{f_{t}}\right) \rho_{\eta_{t}}(x-\varepsilon) d\varepsilon
\label{Prediction}
\end{aligned}
\end{equation}
Here  $\mathcal{Y}_{t}$ denotes the  observations $y_{t}, y_{t-1}, \cdots, y_{0}$.

If not otherwise specified, in the sequel ``prior" refers to the prior density function $\rho_{x_{t+1} \mid \mathcal{Y}_{t}}(x)$ at each time step $t$. We denote by $\mathcal{P}$ the space of all probability density functions supported on $\mathbb{R}$. Let $\mathcal{P}_{2n}$ be the subset of all $\rho\in\mathcal{P}$ which have at least $2n$ finite moments (in addition to $\sigma_0$, which of course is 1). 

We note that the measurement update \eqref{Update2} takes the form of an analytic function with all $\rho_{\epsilon_{t}}$ and $\rho_{x_{t} \mid \mathcal{Y}_{t-1}}$ being non-Gaussian in general. However to obtain an explicit form of prior $\rho_{x_{t+1} \mid \mathcal{Y}_{t}}(x)$ in \eqref{Prediction} when the densities are not Gaussian, is not always a feasible task. We proposed to use the power moments for approximating the prior density in \cite{wu2023non}. For $\rho \in \mathcal{P}_{2n}$, the power moments are calculated by
\begin{equation}
\begin{aligned}
     \sigma_{k, t}
    = & \,\mathbb{E}\left ( x_{t+1}^{k}|\mathcal{Y}_{t} \right )\\
    = & \sum_{j = 0}^{k}\left ( \begin{matrix}
k \\
j
\end{matrix} \right ) f_{t}^{j}\mathbb{E}\left ( x_{t}^{j} |\mathcal{Y}_{t} \right )\mathbb{E}\left ( \eta_{t}^{k-j} \right ),
\end{aligned}
\label{MomentUpdate}
\end{equation}
for $k = 1, \cdots, 2n$ \cite{wu2023non}. The corresponding power moment sequence is denote as $\sigma_{t} = \left(\sigma_{1, t}, \cdots, \sigma_{2n, t} \right)$.

Being linear integral operators, the power moments capture the macroscopic properties of the prior density functions. The success of the power moments in our previous work \cite{wu2023non} naturally leads us to think about using other integral operators to characterize the prior density. Except for the power moments, other statistics have been used in previous research to improve the estimation performance. For example in \cite{2002Identifiability}, covariance lags (power moments) and cepstral coefficients (logarithmic moments) are both used to approximate the spectral density. In this paper, we adopt a similar idea for approximating $\rho_{x_{t+1} \mid \mathcal{Y}_{t}}$ by using both the power and the logarithmic moments.

However, we note that it is not feasible for us to directly use the logarithmic moments in the form of $\int_{\mathbb{R}} x^{k} \log \rho(x) dx$, which are always infinite, even for densities $\rho(x)$ with exponential decaying rate. Take $\rho(x) = \mathcal{N}(0, 1)$ for example. Then
$$
\begin{aligned}
    & \int_{\mathbb{R}} x^{k} \log \left( \frac{1}{\sqrt{2\pi}}\exp\left( -\frac{x^{2}}{2} \right) \right)dx\\
    = & \int_{\mathbb{R}} -\frac{1}{2}\log(2\pi) x^{k} - \frac{1}{2}x^{k+2}dx\\
    = & -\infty, \quad \forall k \in \mathbb{N}_{0}.
\end{aligned}
$$

Therefore, we propose generalized logarithmic moments, for which the first $2n+1$ terms exist and are finite. Thus the generalized logarithmic moments are here defined as
\begin{equation}
\begin{aligned}
    \xi_{k, t} = & \int_{\mathbb{R}} x^{k} \theta(x) \log \rho_{x_{t+1} \mid \mathcal{Y}_{t}}(x) dx\\
    = & \int_{\mathbb{R}} x^{k} \theta(x) \log \int_{\mathbb{R}} \rho_{x_{t} \mid \mathcal{Y}_{t}}\left(\frac{\varepsilon}{f_{t}}\right) \rho_{\eta_{t}}(x-\varepsilon) d\varepsilon dx
\end{aligned}
\label{pik}
\end{equation}
for $k = 1, \cdots, 2n$. They are called ``generalized" because a reference density $\theta(x)$ needs to be determined before calculating them. The corresponding generalized logarithmic moment sequence is denote as $\xi_{t} = \left(\xi_{1, t}, \cdots, \xi_{2n, t} \right)$. We denote by $\mathcal{P}^{\log}_{2n}$ the subset of all $\theta \in \mathcal{P}$ which have finite generalized logarithmic moments to at least order $2n$, provided with $\rho \in \mathcal{P}_{2n}$. Here $\theta \in \mathcal{P}^{\log}_{2n}$ is a reference density function, of which the choice is not very limited indeed. The probability densities with exponential decaying rate, e.g. the exponential families, fall within the subset $\mathcal{P}^{\log}_{2n}$.

Now that the power and generalized logarithmic moments are defined, we give the following definition to characterize the equivalence of two densities in the sense of the two types of moments.

\medskip

\begin{definition} A probability density function, which has the first $2n$ power and generalized logarithmic moment terms identically the same as $\rho$ (with $\theta(x)$ given prior), is called an order-$2n$ P$\&$L density surrogate of $\rho$ and denoted by $\rho^{2n}$.
\end{definition}

We denote by $\hat{\rho}$ the prediction of density $\rho$ and propose to substitute the intractable prior density $\rho_{x_{t+1} \mid \mathcal{Y}_{t}}(x)$ with the proposed density surrogate. Each iteration of Bayesian filtering with the density surrogate is given in Algorithm \ref{Algo1}. At present we assume that it is feasible for us to obtain such a density surrogate given the power and generalized logarithmic moments and first investigate the error propagation through the whole filtering process with the density surrogate, which is one of the most important problems in designing a filter. Since the prior estimation is done at each time step $t$, which means that the approximation errors of the each previous iteration may cause a cumulative one on the current estimation. It distinguishes the problem we treat from conventional density estimation problems. 

\begin{algorithm}[htb] 
\caption{Bayesian filtering with density surrogate at time $t$.} 
\begin{algorithmic}[1] 
\REQUIRE ~~\\ 
System parameters: $f_{t}, h_{t}$;\\
Non-Gaussian densities: $\eta_{t}, \epsilon_{t}$;\\
Prediction at time $t-1$: $\rho_{x_{0}}(x) \text{ or } \hat{\rho}_{x_{t} \mid \mathcal{Y}_{t-1}}(x)$;
\ENSURE ~~\\ 
Prediction at time $t$: $\hat{\rho}_{x_{t+1} \mid \mathcal{Y}_{t}}(x)$;
\STATE Calculate $\hat{\rho}_{x_{t}|\mathcal{Y}_{t}}$ by (\ref{Update1}) or (\ref{Update2}); 
\STATE Calculate $\sigma_{t}$ by (\ref{MomentUpdate});
\STATE Calculate $\xi_{t}$ by (\ref{pik});
\STATE Determine the order-$2n$ P$\&$L density surrogate $\rho^{2n}_{x_{t+1} \mid \mathcal{Y}_{t}}$, for which the truncated power moment sequence is $\sigma_{t}$ and the truncated generalized logarithmic moment sequence is $\xi_{t}$. The prior density estimate at time $t+1$ is then chosen as the P$\&$L density surrogate, i.e., $\hat\rho_{x_{t+1} \mid \mathcal{Y}_{t}} = \rho^{2n}_{x_{t+1} \mid \mathcal{Y}_{t}}$.
\label{AlgoSurrogate}
\end{algorithmic}
\label{Algo1}
\end{algorithm}

We will first review the error propagation of the first $2n$ terms of the power moments in \cite{wu2023non} and then analyze those of the first $2n$ terms of the generalized logarithmic moments. Since the approximation errors caused by the time updates could have cumulative effects on the measurement updates, we analyze the first $2n$ moment terms of not only $\hat{\rho}_{x_{t+1} \mid \mathcal{Y}_{t}}$ but also $\hat{\rho}_{x_{t} \mid \mathcal{Y}_{t}}$.

\medskip

\begin{theorem} 
Suppose $\hat{\rho}_{x_1|\mathcal{Y}_0}$ is a P$\&$L surrogate for $\rho_{x_1|\mathcal{Y}_0}$, and let $\hat{\rho}_{x_t|\mathcal{Y}_t}$ and $\hat{\rho}_{x_{t+1}|\mathcal{Y}_t}$ be obtained from Algorithm 1 for $t = 2, 3, \cdots$. Then the power moments and the generalized logarithmic moments of $\hat{\rho}_{x_t|\mathcal{Y}_t}$ and $\hat{\rho}_{x_{t+1}|\mathcal{Y}_t}$ are asymptotically unbiased from those of $\rho_{x_t|\mathcal{Y}_t}$ and $\rho_{x_{t+1}|\mathcal{Y}_t}$, respectively, and are approximately identical to them for a sufficiently large $n$, given that all power moments and generalized logarithmic moments of $x_{t}$ and the corresponding $\hat{x}_{t}$ exist and are finite.
\label{MomentError}
\end{theorem}

A complete proof of Theorem \ref{MomentError} is given in Appendix \ref{Appendix_asymptotic}. Theorem \ref{MomentError} reveals the fact that the first $2n$ generalized logarithmic moment terms of the estimated prior densities with the density surrogate are approximately identical to the true ones through the whole filtering process. Together with \eqref{MomentApprox}, we have that $\hat \rho_{x_{t+1} \mid \mathcal{Y}_{t}}$ and $\hat \rho_{x_{t} \mid \mathcal{Y}_{t}}$ are approximately order-$2n$ P$\&$L density surrogates of $\rho_{x_{t+1} \mid \mathcal{Y}_{t}}$ and $\rho_{x_{t} \mid \mathcal{Y}_{t}}$. It reveals the fact that approximation using both moments does not introduce significant cumulative errors to the first $2n$ moment terms of the estimated pdfs, with $n$ chosen as a relatively large integer.

The problem is now constructing an order-$2n$ P$\&$L density surrogate. Since the domain of $\rho$ is $\mathbb{R}$, the problem becomes a Hamburger moment problem \cite{schmudgen2017moment} with the constraints of additive generalized logarithmic moments. In the next section, we will give a representation of this specific moment problem and propose a solution to it.

\section{A parametrization of the density surrogate using power and generalized logarithmic moments}
\label{ParameterizationOf}

In this section, we give a formal definition of the approximation problem of the prior density and prove the existence of a solution to this problem given the power moments and the generalized logarithmic moments. 

\medskip

\begin{definition}
 A sequence 
 $$
 (\sigma_1, \ldots, \sigma_{2n}, \xi_1, \ldots, \xi_{2n})
 $$
 is a feasible $2n$ P$\&$L sequence, if there is a random variable $X$ with a probability density $\rho(x)$ defined on $\mathbb{R}$, whose moments are given by \eqref{MomentUpdate} and \eqref{pik}, that is,
\begin{equation*}
    \sigma_{k} = \mathbb{E}\{X^{k} \} =\int_\mathbb{R}x^k\rho(x)dx, \quad k = 1,\dots, 2n,
\end{equation*}
and
\begin{equation*}
    \xi_{k} = \mathbb{E}^{\log}\{X^{k} \} = \int_{\mathbb{R}} x^{k} \theta(x) \log \rho(x) dx, k = 1, \cdots, 2n.
\end{equation*}

Moreover, we assume that $\sigma_{0} =1$ and $\xi_{0} = 0$. Any such random variable $X$ is said to have a $(\sigma, \xi)$-feasible distribution. We denote the random variable as $X \sim(\sigma, \xi)$. 
\label{DefinitionMoment}
\end{definition}

Next we prove the existence of a solution to the moment problem defined in Definition \ref{DefinitionMoment}. We first paraphrase Exercise 13.12 in \cite{rudin_2015}. Let $f$ be a real-valued measurable function defined on $\mathbb{R}$. Then there exists a sequence of polynomials $P_{n}$ such that
$$
\lim_{n \rightarrow +\infty}P_{n}(x)=f(x)$$ 
almost everywhere.

We note that the true $\rho_{x_{t+1} \mid \mathcal{Y}_{t}}(x)$ is trivially a solution to the moment problem in Definition \ref{DefinitionMoment}. However we require an analytic function which satisfy the moment constraints. In our problem setting, $\rho_{x_{t+1}|\mathcal{Y}_{t}}$ is Lebesgue measurable. Therefore there exists a $\lim_{n \rightarrow +\infty}P_{n}(x)$ which is equal to $\rho_{x_{t+1} \mid \mathcal{Y}_{t}}(x)$ almost everywhere, i.e., it is a solution to the moment problem above.

However this solution does not exactly satisfy the requirement of a state estimate of the Bayesian filter, since there are possibly infinitely many parameters in the solution, which makes it infeasible to propagate the solution in the filtering process. Parametrization is then the most significant problem, which aims to use finitely many parameters to characterize the density. 

Meanwhile, we are provided with two truncated power and generalized logarithmic moment sequences rather than two full ones, which means that there might be infinitely many feasible solutions to this problem. In the following part of this section, we propose to choose proper constraints to parameterize the density surrogate that satisfies the moment conditions. We still emphasize here that the parametrization is not unique. Different constraints will yield different parametrizations.

In the following part of this section, we propose to parameterize the density surrogate, i.e., to derive a unique solution to the moment problem of $\rho_{x_{t+1} \mid \mathcal{Y}_{t}}$. For simplicity, we omit the subscript $t$ in all the terms in the following part of this section.

Since there are infinitely many feasible solutions to the moment problem, a criterion to determine a unique solution is necessary. Following \cite{georgiou2003kullback, wu2023non}, we consider the Kullback-Leibler (KL) distance
\begin{equation}
\label{KL}
\mathbb{KL}(\theta\|\rho)=\int_\mathbb{R} \theta(x) \log \frac{\theta(x)}{\rho(x)} dx
\end{equation}
to measure the difference between $\theta$ and $\rho$, which is a widely used pseudo-measure in density estimation tasks \cite{georgiou2003kullback, hall1987kullback, li1999mixture, vapnik1999nature}. Although it is not symmetric, which makes it not a real metric, the KL distance is jointly convex. We formulate the primal problem as minimizing
\begin{equation}
\int_{\mathbb{R}} \theta(x) \log \frac{\theta(x)}{\rho(x)} d x
\label{Primal}
\end{equation}
with respect to $\rho(x)$ and subject to
\begin{equation}
\int_{\mathbb{R}} x^k \rho d x=\sigma_k, \quad k=1, \ldots, 2 n
\label{Constraint1}
\end{equation}
and
\begin{equation}
\int_{\mathbb{R}} x^k \theta \log \rho d x=\xi_k, \quad k=1, \ldots, 2 n.
\label{Constraint2}
\end{equation}

Here $\theta$ is a prespecified density function which we want the estimate of the prior density $\rho_{x_{t+1} \mid \mathcal{Y}_{t}}(x)$ to be close to. As we have mentioned, there are infinitely many solutions to the truncated moment problem, among which some do not have satisfactory properties (e.g. smoothness) or have undesired massive modes (peaks). By minimizing the Kullback-Leibler distance between $\theta$ and the density estimate, with $\theta$ given as an analytic function, it is possible to  obtain an analytic density estimate of the prior, as will be proved in the following part of this section. We are then able to propagate the prior density throughout the filtering process. By minimizing \eqref{Primal} subject to \eqref{Constraint1} and \eqref{Constraint2}, a parametrization based on both the power moments and the generalized logarithmic moments is proposed in the following theorem.

\medskip

\begin{theorem}
\label{theorem32}
Denoting the Lagrange multipliers as
$$p=(p_1,p_2,\dots,p_{2n}), \quad q=(q_1,\dots,q_{2n}),
$$
set
\begin{displaymath}
P(x)=1+p_1x+ p_2x^2 +\dots +p_{2n}x^{2n},
\end{displaymath}
\begin{displaymath}
Q(x)=q_0+q_1x+ q_2x^2 +\dots +q_{2n}x^{2n}.
\end{displaymath}

Given any $\theta \in \mathcal{P}^{\log}_{2n}$ and any $\sigma$ which satisfies 
\begin{equation}
\begin{bmatrix}
1 & \sigma_{1} & \cdots & \sigma_{n}\\ 
\sigma_{1} & \sigma_{2} & \cdots & \sigma_{n+1}\\ 
\vdots &  & \ddots & \\ 
\sigma_{n} & \sigma_{n+1} &  & \sigma_{2n}
\end{bmatrix} \succ 0,
\label{PosDef}
\end{equation}
minimizing \eqref{Primal}
subject to \eqref{Constraint1} and \eqref{Constraint2} yields a unique solution $\rho \in \mathcal{P}_{2n}$ of the form 
\begin{equation}
    \hat{\rho}=\frac{\hat{P}(x)}{\hat{Q}(x)}\theta,
\label{MiniForm}
\end{equation}
where $\hat{P}(x)$ and $\hat{Q}(x)$ are coprime and $\hat{p}, \hat{q}$ corresponding to $\hat{P}(x)$ and $\hat{Q}(x)$ are the unique solutions to the problem of minimizing
\begin{equation}
\mathbb{J}(P, Q)=\sigma'q  - \xi'p +\int_{\mathbb{R}} P\theta\log\frac{P\theta}{Q} dx - \int_{\mathbb{R}} P\theta dx.
\label{LossFuncpq}
\end{equation}
over all $P(x), Q(x) > 0$.
\end{theorem}

\begin{proof}

The Lagrangian of the primal problem \eqref{Primal} with constraints \eqref{Constraint1} and \eqref{Constraint2} is
\begin{displaymath}
\begin{aligned}
L(\rho,p,q) = & \int_{\mathbb{R}}\theta\log\frac{\theta}{\rho}dx +\sum_{k=0}^{2n}q_k\left(\int_{\mathbb{R}} x^k\rho dx -\sigma_k\right)\\
- & \sum_{k=1}^{2n}p_k\left(\int_{\mathbb{R}} x^k\theta\log\rho dx -\xi_k\right).
\end{aligned}
\end{displaymath}

Then we have 
\begin{align*}
   L(\rho,p,q) = & \int_{\mathbb{R}}\theta\log\frac{\theta}{\rho}dx -\int_{\mathbb{R}} (P-1)\theta\log\rho dx\\
   + & \xi' p + \int_{\mathbb{R}} Q\rho dx -\sigma'q \\
   = & \int_{\mathbb{R}}\theta\log\theta dx -\int_{\mathbb{R}} P\theta\log\rho dx\\
   + & \int_{\mathbb{R}} Q\rho dx +\xi'p -\sigma'q  
\end{align*}
with the directional derivative
\begin{displaymath}
\delta L(\rho,p,q;\delta\rho) = \int_{\mathbb{R}}\delta\rho\left(Q-\frac{P\theta}{\rho}\right)dx
\end{displaymath}
which yields the stationary point 
\begin{displaymath}
\hat{\rho} =\frac{P\theta}{Q}.
\end{displaymath}
Then
\begin{displaymath}
L(\hat{\rho},p,q) = \int_{\mathbb{R}}\theta\log\theta dx - \mathbb{J}(P, Q),
\end{displaymath}
where
\begin{displaymath}
\mathbb{J}(P, Q)=\sigma'q  - \xi'p +\int_{\mathbb{R}} P\theta\log\frac{P\theta}{Q} dx - \int_{\mathbb{R}} P\theta dx.
\end{displaymath}
In particular,
\begin{align*}
\frac{\partial\mathbb{J}}{\partial p_k}    &=  -\xi_k +\int_{\mathbb{R}} x^k\theta\log\frac{P\theta}{Q}dx \\
  \frac{\partial\mathbb{J}}{\partial q_k}  &  =\sigma_k - \int_{\mathbb{R}} x^k\frac{P\theta}{Q}dx.
\end{align*}
\end{proof}

\begin{remark}
We have choosen the constant term of $P(x)$ as $1$ to yield a simpler form of $\hat{\rho}$ in \eqref{MiniForm}. However, it can be any real number. By specifying
\begin{equation}
\int_{\mathbb{R}} \theta\log\frac{P\theta}{Q}dx = 0
\label{xi0}
\end{equation}
and
\begin{equation}
    \int_{\mathbb{R}} \frac{P\theta}{Q}dx =1,
\label{sigma0}
\end{equation}
it is determined together with $q_{0}$.
\end{remark}

But here we note that our parametrization has a rational form. Therefore, by dividing both the numerator and denominator by $p_{0}$, the constant term of the numerator becomes $1$, and the $q_{0}/p_{0}$ becomes the new $q_{0}$ in Theorem \ref{theorem32}. 

We note that to obtain a unique solution to the problem to minimize \eqref{LossFuncpq} by a gradient-based method, it remains to prove that the map from $(p, q)$ to $(\xi, \sigma)$ is a diffeomorphism. In the following section, we will complete the proof of Theorem \ref{theorem32} by proving precisely this.

\section{The diffeomorphic map}
\label{AGlobal}

In this section, we prove that the map $(P, Q) \mapsto (\sigma, \xi)$ is diffeomorphic, building upon some of the ideas presented in \cite{2002Identifiability}.

We begin by noting that $\sigma_0$ is always equal to one. We also consider $q_0$ as a normalizing factor to ensure that $\int_{\mathbb{R}} \rho(x) dx = 1$, which is thus determined when $(p_1, p_2, \dots, p_{2n})$ and $(q_1, q_2, \dots, q_{2n})$ are known. Therefore, denoting by $\mathcal{S}_{2n}$ as the class of positive polynomials of order $2n$ with the term of order zero being a constant, we have $P, Q \in \mathcal{S}_{2n}$. Given a specified density function $\theta(x) \in \mathcal{P}^{\log}_{2n}$, we can represent the rational density function by $(P, Q) \in \mathcal{M}_{2n}$, where $\mathcal{M}_{2n} = \mathcal{S}_{2n} \times \mathcal{S}_{2n}$. Thus $\mathcal{M}_{2n}$ becomes a smooth, connected, real manifold of dimension $4n$ which is diffeomorphic to $\mathbb{R}^{4n}$.

Next, we define some additional spaces for analysis. We denote by $\mathcal{M}_{2n}^{*}$ the (dense) open subspace of $\mathcal{M}_{2n}$ consisting of pairs $(P, Q)$ of coprime polynomials. For $P \in \mathcal{S}_{2n}$, we define $\mathcal{M}_{2n}(P)$ as the space of all points in $\mathcal{M}_{2n}$ with the polynomial $P$ fixed. Similarly defining $\mathcal{M}_{2n}(Q)$,  $\mathcal{M}_{2n}(P)$ and $\mathcal{M}_{2n}(Q)$ become real, smooth, connected $2n$-manifolds that are diffeomorphic to $\mathcal{S}_{2n}$ and thus to $\mathbb{R}^{2n}$. Furthermore, the tangent vectors to $\mathcal{M}_{2n}$ at $(P, Q)$ can be represented as perturbations $(P + \epsilon u, Q + \epsilon v)$, where $u$ and $v$ are polynomials of degree less than or equal to $2n - 1$. Denoting the real vector space of polynomials of degree less than or equal to $d$ by $V_d$, the tangent space to $\mathcal{M}_{2n}$ at a point $(P, Q)$ is canonically isomorphic to $V_{2n-1} \times V_{2n-1}$. Additionally, the tangent space to the submanifold $\mathcal{M}_{2n}(P)$ at a point $(P, Q)$ is given by
\begin{equation*}
    T_{(P, Q)} \mathcal{M}_{2n}(P)=\left\{(u, v) \in V_{2n-1} \times V_{2n-1} \mid u=0\right\}
\end{equation*}

Similarly, the tangent space to $\mathcal{M}_{2n}(Q)$ is given by
\begin{equation*}
    T_{(P, Q)} \mathcal{M}_{2n}(Q)=\left\{(u, v) \in V_{2n-1} \times V_{2n-1} \mid v=0\right\}    
\end{equation*}

The $2n$-manifolds $\{\mathcal{M}_{2n}(P) \mid P \in \mathcal{S}_{2n}\}$ form the leaves of a foliation of $\mathcal{M}_{2n}$, as do the $2n$-manifolds $\{\mathcal{M}_{2n}(Q) \mid Q \in \mathcal{S}_{2n}\}$. Furthermore, these two foliations are complementary in the sense that if a leaf of one intersects a leaf of the other, the tangent spaces intersect only at $(0,0)$. This transversality property is equivalent to the fact that the polynomials $(P, Q)$ form a local system of coordinates.

From a geometric perspective, this property implies that $(P, Q)$ are smooth coordinates on $\mathcal{M}_{2n}$. We will use this to demonstrate that $(P, Q)$ also form bona-fide coordinate systems. Let $g: \mathcal{M}_{2n} \rightarrow \mathbb{R}^{2n}$ be the map that sends $(P, Q)$ to $\xi$, where the components of $\xi$ are calculated using equation \eqref{pik}. We denote $\mathcal{C}_{2n}:=g(\mathcal{M}_{2n})$. Additionally, for each $\xi \in \mathcal{C}_{2n}$, we define the subset $\mathcal{M}_{2n}(\xi)=g^{-1}(\xi)$.

We aim to show that $\mathcal{M}_{2n}(\xi)$ is a smooth submanifold of dimension $2n$. To achieve this, we need to compute the Jacobian matrix of $g$ evaluated at tangent vectors to a point $(P, Q) \in \mathcal{M}_{2n}$. If the Jacobian matrix of $g$ is full rank at every point $(P, Q) \in \mathcal{M}_{2n}$, meaning that the directional derivative exists in every direction at each point, then $\mathcal{M}_{2n}(\xi)$ is proved to be smooth.

We recall that the tangent vectors to $\mathcal{M}_{2n}$ at $\left ( P, Q \right )$  can be expressed as a perturbation $\left ( P + \epsilon u , Q + \epsilon v \right )$, where $u, v$ are polynomials of degree less than or equal to $2n-1$. For each component ($k = 1, \cdots, 2n$)
\begin{equation}
    g_{k}\left ( P, Q \right ) = \int_{\mathbb{R}}x^{k}\theta(x)\log \left( \frac{P(x)}{Q(x)} \theta(x)\right )dx
\end{equation}
of $g$, we construct the directional derivative as follows: 
\begin{equation}
\begin{aligned}
    & D_{\left ( u, v \right )} g_{k}\left ( P, Q \right )\\
    & = \lim_{\epsilon \rightarrow 0} \frac{1}{\epsilon} \left [ g_{k}\left ( P + \epsilon u, Q + \epsilon v \right ) - g_{k}\left ( P, Q \right )\right ]\\
    & = \int_{\mathbb{R}}\left ( \frac{u}{P} - \frac{v}{Q}\right )\theta x^{k}dx
\end{aligned}
\label{derivative}
\end{equation}
in the direction $(u,v) \in V_{2n-1} \times V_{2n-1}$.

Next, we define the linear map $G: V_{2n-1} \mapsto \mathbb{R}^{2n}$ as follows:

\begin{equation}
    G_{\psi}u = \int_{\mathbb{R}} \frac{u}{\psi}\theta\begin{bmatrix}
x\\ 
x^{2}\\ 
\vdots \\ 
x^{2n}
\end{bmatrix}dx.
\label{G}
\end{equation}

Then the kernel of the Jacobian of $g$ at $\left ( P, Q \right )$ is given by
\begin{equation}
    \ker \text{Jac}(g)|_{(u,v)} = \left \{ (P, Q)|G_{P}u = G_{Q}v \right \}
\label{kernel}
\end{equation}


\medskip

\begin{lemma}
The linear map $G_{\psi}$ is a bijection.
\label{lemma51}
\end{lemma}
\begin{proof}
First, consider the case when $G_{\psi}u = 0$. This implies that
\begin{equation}
\int_{\mathbb{R}} \frac{u}{\psi} \theta x^{k}dx = 0
\label{linop}
\end{equation}
for $k = 1, \cdots, 2n$. 

From \eqref{xi0}, we have that $g_{0}\left( \psi, Q \right) \equiv 0$ for any $(\psi, Q) \in \mathcal{S}_{2n}$, which means that directional derivative along any direction is equal to zero. We take the directional derivative along $(u, 0)$, and we have
\begin{equation}
\begin{aligned}
    & D_{\left ( u, 0 \right )} g_{0}\left ( \psi, Q \right )\\
    & = \lim_{\epsilon \rightarrow 0} \frac{1}{\epsilon} \left [ g_{k}\left ( \psi + \epsilon u, Q \right ) - g_{k}\left ( \psi, Q \right )\right ]\\
    & = \int_{\mathbb{R}}\frac{u}{\psi}\theta dx = 0.
\end{aligned}
\label{g0deriv}
\end{equation}
Since $u \in V_{2n-1}$, we write
$$
u(x) = \sum_{i=0}^{2n-1} u_{i}x^{i}, u_{i} \in \mathbb{R}.
$$

By \eqref{linop} and \eqref{g0deriv}, we shall write
$$
\sum_{i = 0}^{2n-1} u_{i} \int_{\mathbb{R}} \frac{u}{\psi} \theta x^{i}dx = \int_{\mathbb{R}} \frac{u^{2}}{\psi} \theta dx = 0.
$$
Since $\theta, \psi$ are both positive, $u(x)$ needs to be zero.  Therefore, we have established the injectivity of $G_{\psi}$. Furthermore, since the range and domain of $G_{\psi}$ have the same dimension, namely $2n$, the map is also surjective. Consequently, we can conclude that $G_{\psi}$ is a bijection.
\end{proof}

\begin{proposition} For each  $\xi \in \mathcal{C}_{2n}$, the space $\mathcal{M}_{2n}(\xi)$ is a smooth $2n$-manifold. The tangent space $T_{(P, Q)} \mathcal{M}_{2n}(\xi)$ at $(P, Q)$ consists of precisely all $(u, v) \in V_{2n-1} \times V_{2n-1}$ such that
\begin{equation}
    \int_{\mathbb{R}} \frac{u}{P} \theta x^{k} dx= \int_{\mathbb{R}} \frac{v}{Q} \theta x^{k} dx
\label{uvPQ}
\end{equation}
for $k=0,1, \ldots, 2n$.
\label{proposition52}
\end{proposition}
\begin{proof}
The tangent vectors of $\mathcal{M}_{2n}(\xi)$ at $(P, Q)$ correspond to the vectors in the null space of the Jacobian of $g$ at $(P, Q)$, as indicated by equation (\ref{derivative}). Consequently, by utilizing equation \eqref{kernel}, we establish that \eqref{uvPQ} holds for $k = 1, 2, \ldots, 2n$. Additionally, according to Lemma \ref{lemma51}, we can conclude that \eqref{linop} holds for $k = 0$. Therefore, \eqref{uvPQ} also holds for $k = 0$. Furthermore, based on \eqref{kernel} and Lemma \ref{lemma51}, the tangent space has a dimension of $2n$. Consequently, the Jacobian matrix $\left.\operatorname{Jac}(g)\right|{(P, Q)}$ has full rank, and the remaining part of the claim follows from the implicit function theorem.

Since the rank of $\left.\operatorname{Jac}(g)\right|{(P, Q)}$ is consistently $2n$, the connected components of the submanifolds $\mathcal{M}_{2n}(\xi)$ constitute the leaves of a foliation of $\mathcal{M}_{2n}$. However, we still need to demonstrate that the submanifolds $\mathcal{M}_{2n}(\xi)$ themselves are connected. The detailed proof for this is provided in Appendix \ref{Append_sigma}. Consequently, we can state the following proposition.
\end{proof}

\begin{proposition} The $2n$-manifolds $\left\{\mathcal{M}_{2n}(\xi) \mid \xi \in \mathcal{C}_{2n}\right\}$ are connected, hence forming the leaves of a foliation of $\mathcal{M}_{2n}$.
\end{proposition}

From the results proved so far, we conclude the following corollary.

\medskip

\begin{corollary}
The foliations, $\left\{\mathcal{M}_{2n}({Q}) \mid Q \in \mathcal{S}_{2n}\right\}$ and $\left\{\mathcal{M}_{2n}({\xi}) \mid {\xi} \in \mathcal{C}_{2n} \right\}$, are
complementary, i.e., any intersecting pair of leaves, with one leaf from each foliation, intersects transversely. And each intersecting pair of leaves intersects in at most one point.
\end{corollary}

\begin{proof}
Setting $u=0$ in (\ref{kernel}), we obtain $G_{Q} v=0 .$ Hence, by Lemma~\ref{lemma51}, $v=0$ so that the foliations are transverse. If a leaf $\mathcal{P}_{2n}(P)$ intersects a leaf $\mathcal{M}_{2n}(\xi)$ at a point $(P, Q)$, then the corresponding $P$ is known. Then according to Appendix \ref{Append_xi}, a unique $\xi$ is determined.
\end{proof}

A similar statement for the foliation $\left\{\mathcal{M}_{2n}(Q) \mid Q \in \mathcal{S}_{2n}\right\}$ can be proved by the mirror image of this proof and will be omitted.

Next, let $h: \mathcal{M}_{2n} \rightarrow \mathbb{R}^{2n}$ be the map which sends $(P, Q)$ to $\sigma$, the components of which are calculated by \eqref{MomentUpdate}, and let $\mathcal{R}_{2n} = h\left ( \mathcal{M}_{2n} \right )$. Each $\sigma$ satisfies \eqref{PosDef}, guaranteeing the existence of a solution to the moment problem.

Now for each $\sigma \in \mathcal{R}_{2n}$, we aim to demonstrate that the set
\begin{equation}
    \mathcal{M}_{2n}\left ( \sigma \right ) = h^{-1}\left ( \sigma \right )
\end{equation}
forms a smooth manifold of dimension $2n$. The tangent vectors to $\mathcal{M}_{2n}$ at $\left ( P, Q \right )$ can be represented as a perturbation $\left ( P + \epsilon u , Q + \epsilon v \right )$, where $u, v$ are polynomials of degree less than or equal to $2n-1$. For each component
\begin{equation}
    h_{k}\left ( P, Q \right ) = \int_{\mathbb{R}}x^{k}\theta(x)\frac{P(x)}{Q(x)}dx, \quad k = 1, \cdots, 2n,
\end{equation}
the directional derivative of $h$ at $(P, Q) \in$ $\mathcal{M}_{2n}$ in the direction $(u, v) \in V_{2n-1} \times V_{2n-1}$ is
\begin{equation}
\begin{aligned}
    & D_{\left ( u, v \right )} h_{k}\left ( P, Q \right )\\
    & = \lim_{\epsilon \rightarrow 0} \frac{1}{\epsilon} \left [ h_{k}\left ( P + \epsilon u, Q + \epsilon v \right ) - h_{k}\left ( P, Q \right )\right ]\\
    & = \int_{\mathbb{R}}\left ( \frac{u Q}{Q^2} - \frac{v P}{Q^2} \right )\theta x^{k}dx
\end{aligned}
\label{derivativeMu}
\end{equation}

Similar to (\ref{G}), we define the linear map $H_{\psi}: V_{2n-1} \rightarrow \mathbb{R}^{2n}$ by
\begin{equation}
    H_{\psi}u = \int_{\mathbb{R}} \frac{u \psi}{Q^2}\theta \begin{bmatrix}
x\\ 
x^{2}\\ 
\vdots \\ 
x^{2n}
\end{bmatrix}dx
\end{equation}
and the kernel of the Jacobian of $h$ at $(P, Q)$ is given by
\begin{equation}
    \ker \text{Jac}(h)|_{(P, Q)} = \left \{ (u, v)|H_{Q}u = H_{P}v \right \}
\end{equation}
for $k=0,1, \ldots, 2n$.

\medskip

\begin{proposition}
For each $\sigma \in \mathcal{R}_{2n}$, the subspace $\mathcal{M}_{2n}(\sigma)$ is a smooth and connected $2n$-manifold. The tangent space $T_{(P, Q)} \mathcal{M}_{2n}(\sigma)$ consists of pairs $(u, v) \in V_{2n-1} \times V_{2n-1}$ satisfying
\begin{equation}
\int_{\mathbb{R}} \frac{u Q}{Q^2} \theta x^{k} dx= \int_{\mathbb{R}} \frac{v  P}{Q^2} \theta x^{k} dx
\label{uvQP}
\end{equation}
for $k = 1, \ldots, 2n$. Furthermore, the $2n$-manifolds $\mathcal{M}_{2n}(\sigma)$ constitute the leaves of a foliation of $\mathcal{M}_{2n}$.
\label{proposition55}
\end{proposition}

\begin{proof}
Let us begin by demonstrating that the linear map $H_{\psi}$ is a bijection. Suppose $H_{\psi}u = 0$. This implies
\begin{equation}
H_{\psi}u = \int_{\mathbb{R}} \frac{u \psi}{Q^2} \theta x^{k}dx = 0
\label{linopH}
\end{equation}
for $k = 0, \ldots, 2n$. From \eqref{sigma0}, we have $h_{0}\left( \psi, Q \right) = 1$ for any $(\psi, Q) \in \mathcal{S}_{2n}$. Therefore, the directional derivative along any direction is equal to zero. We take the directional derivative along $(u, 0)$, and we have
\begin{equation}
\begin{aligned}
    & D_{\left ( u, 0 \right )} h_{0}\left ( \psi, Q \right )\\
    & = \lim_{\epsilon \rightarrow 0} \frac{1}{\epsilon} \left [ h_{k}\left ( \psi + \epsilon u, Q \right ) - h_{k}\left ( \psi, Q \right )\right ]\\
    & = \int_{\mathbb{R}} \frac{u  \psi}{Q^2} \theta dx = 0.
\end{aligned}
\label{h0deriv}
\end{equation}
Since $u \in V_{2n-1}$, we write
$$
u(x) = \sum_{i=0}^{2n-1} u_{i}x^{i}, u_{i} \in \mathbb{R}.
$$

By \eqref{linopH} and \eqref{h0deriv}, we shall write
$$
\sum_{i = 0}^{2n-1} u_{i} \int_{\mathbb{R}} \frac{u  \psi}{Q^2} \theta x^{i}dx = \int_{\mathbb{R}} \frac{u^{2}  \psi}{Q^{2}} \theta dx = 0.
$$
Since $\theta, \psi$ are both positive, we conclude that $u = 0$. Thus, $H_{\psi}$ is injective. Moreover, since the range and domain of $H_{\psi}$ have the same dimension, namely $2n$, the map is also surjective. In conclusion, $H_{\psi}$ is a bijection.
Similar to Proposition \ref{proposition52}, we can establish that the rank of $\left.\operatorname{Jac}(h)\right|{(P, Q)}$ is full. As the rank of $\left.\operatorname{Jac}(h)\right|{(P, Q)}$ is thus consistently $2n$, the connected components of the submanifolds $\mathcal{M}_{2n}(\sigma)$ form the leaves of a foliation of $\mathcal{M}_{2n}$. To complete the argument, it remains to demonstrate that the submanifolds $\mathcal{M}_{2n}(\sigma)$ are themselves connected. The proof is provided in Appendix \ref{Append_xi}.
\end{proof}

\begin{theorem}
For each $(P, Q) \in \mathcal{M}_{2n}(\sigma) \cap \mathcal{M}_{2n}(\xi)$, the dimension of
\begin{equation}
   \mathcal{D} :=T_{(P, Q)} \mathcal{M}_{2n}(\sigma) \cap T_{(P, Q)} \mathcal{M}_{2n}(\xi) 
\end{equation}
equals the degree of the greatest common divisor of the polynomials $P(x)$ and $Q(x)$.
\label{theorem56}
\end{theorem}

\begin{proof} Every $(P, Q) \in \mathcal{D}$ satisfies both (\ref{uvPQ}) and (\ref{uvQP}). By taking appropriate linear combinations of (\ref{uvPQ}) and (\ref{uvQP}), we obtain the following equations
\begin{equation}\label{uP}
    \int_{\mathbb{R}} \frac{u^{2}}{P^2}\theta dx = \int_{\mathbb{R}} \frac{uv}{PQ}\theta dx,
\end{equation}
and
\begin{equation}\label{vQ}
   \int_{\mathbb{R}} \frac{u v}{P Q}\theta dx = \int_{\mathbb{R}} \frac{v^{2}}{Q^2}\theta dx.
\end{equation}

Since $\theta(x)$ is a non-negative density function, we can define
\begin{equation}
    f_{1} := \frac{u}{P}\theta^{\frac{1}{2}} \quad \text{and} \quad f_{2} := \frac{v}{Q}\theta^{\frac{1}{2}},
\end{equation}
which allows us to rewrite \eqref{uP} and \eqref{vQ} as
\begin{equation}
\left\|f_{1}\right\|^{2}=\left\langle f_{1}, f_{2}\right\rangle \quad \text { and } \quad\left\langle f_{1}, f_{2}\right\rangle=\left\|f_{2}\right\|^{2}
\end{equation}
using the inner product and norm of $L^{2}[-\infty, +\infty]$. Applying the parallelogram law, we have
\begin{equation}
    \left\|f_{1}-f_{2}\right\|^{2}=\left\|f_{1}\right\|^{2}+\left\|f_{2}\right\|^{2}-2\left\langle f_{1}, f_{2}\right\rangle=0,    
\end{equation}
which implies $f_{1}=f_{2}$. Consequently,
\begin{equation}
    \frac{u}{v} = \frac{P}{Q}
\end{equation}
which has no solution if $P$ and $Q$ are coprime. However, if $P$ and $Q$ have a greatest common factor of degree $d$, $u(x)$ and $v(x)$ can be polynomials of degree less than or equal to $2n-1$ with an arbitrary common factor of degree $d-1$, thereby defining a vector space of dimension $d$, as stated.
\end{proof}

In conclusion, Theorem \ref{theorem56} establishes the complementary property of the foliations $\left\{\mathcal{M}_{2n}(\sigma) \mid \sigma \in \mathcal{L}_{+} \right\}$ and $\left\{\mathcal{M}_{2n}(\xi) \mid \xi \in \mathcal{L}_{+} \right\}$ at every point $(P, Q) \in \mathcal{M}_{2n}^{*}$, where $P$ and $Q$ are coprime. Consequently, it follows that the kernels of $\left.\operatorname{Jac}(g)\right|{(P, Q)}$ and $\left.\operatorname{Jac}(h)\right|{(P, Q)}$ are complementary at any point $(P, Q)$ in $\mathcal{M}_{2n}^{*}$. Remarkably, the Jacobian of the joint map $(P, Q) \mapsto \left ( \sigma, \xi \right )$ achieves full rank. As a result, the mapping $(P, Q) \mapsto (\sigma, \xi)$ is a diffeomorphism, thereby completing the proof of the following theorem.

\medskip

\begin{theorem}
The power moments $\sigma_{1}, \sigma_{2}, \cdots, \sigma_{2n}$ and the generalized logarithmic moments $\xi_{1}, \xi_{2}, \cdots, \xi_{2n}$ serve as a valid smooth coordinate system within the open subset $\mathcal{M}_{2n}^{*}$ of $\mathcal{M}_{2n}$. This means that the mapping from $\mathcal{M}_{2n}^{*}$ to $\mathbb{R}^{4n}$ with components
$$
\left(\sigma_{1}, \sigma_{2}, \cdots, \sigma_{2n}, \xi_{1}, \xi_{2}, \cdots, \xi_{2n}\right)
$$
has an everywhere invertible Jacobian matrix.
\label{theorem57}
\end{theorem}

Based on the conclusive findings presented in Theorem \ref{theorem57}, we have now concluded the proof for Theorem \ref{theorem32}.

Having obtained all the necessary results from the preceding sections, we can now provide a comprehensive algorithm for non-Gaussian Bayesian filtering utilizing moments. This algorithm, denoted as Algorithm \ref{Algo2}, uses both types of moments and is built upon the foundation of Algorithm \ref{Algo1}.

\begin{algorithm}[htb] 
\caption{Bayesian filtering with density surrogate using power moments at time $t$.} 
\begin{algorithmic}[1] 
\REQUIRE ~~\\ 
System parameters: $f_{t}, h_{t}$;\\
Non-Gaussian densities: $\eta_{t}, \epsilon_{t}$;\\
Prediction at time $t-1$: $\rho_{x_{0}}(x) \text{or } \hat{\rho}_{x_{t} \mid \mathcal{Y}_{t-1}}(x)$;
\ENSURE ~~\\ 
Prediction at time $t$: $\hat{\rho}_{x_{t+1} \mid \mathcal{Y}_{t}}(x)$;
\STATE Calculate $\hat{\rho}_{x_{t}|\mathcal{Y}_{t}}$ by (\ref{Update1}) or (\ref{Update2}); 
\STATE Calculate $\sigma_{t}$ by (\ref{MomentUpdate});
\STATE Calculate $\xi_{t}$ by (\ref{pik});
\STATE Perform optimization, solving \eqref{LossFuncpq} to obtain the order-$2n$ $P \& L$ density surrogate of $\int_{\mathbb{R}} \hat \rho_{x_{t} \mid \mathcal{Y}_{t}}\left(\frac{\varepsilon}{f_{t}}\right) \rho_{\eta_{t}}(x-\varepsilon) d\varepsilon$, which represents the new prediction $\hat{\rho}_{x_{t+1} \mid \mathcal{Y}_{t}}(x)$.
\end{algorithmic}
\label{Algo2}
\end{algorithm}

In \cite{wu2023non}, we demonstrated that power moments, which are linear integral operators, contain abundant information for characterizing density functions and can transform the infinite-dimensional filtering problem into a finite-dimensional and tractable one. In this paper, we prove that other linear integral operators, which capture different types of macroscopic properties of the density to be estimated, provide additional information that can enhance the density estimate. In the next section, we will simulate our proposed non-Gaussian Bayesian filter on mixtures of different types of density functions. We will also compare its performance with the filter proposed in \cite{wu2023non} for each numerical example, demonstrating that the additional information carried by the generalized logarithmic moments improves the density estimation performance. Moreover, we will apply the proposed filter to a robot localization task and compare the performance to several prevailing methods.

\section{Numerical examples}
In this section, we provide numerical examples of our proposed non-Gaussian Bayesian filter that utilizes both the power moments and the generalized logarithmic moments.

We perform two types of numerical simulations for validating the performance of the proposed algorithm. We first simulate distribution approximation tasks. We compare this filter, denoted as DPBM (Density Parametrization using both Power Moments and Generalized Logarithmic Moments), with a Bayesian filter that only uses power moments, which we referred to as DPPM (Density Parametrization using Power Moments) in our previous paper \cite{wu2023non}.

To begin, we need to choose a reference density $\theta(x)$. For light-tailed density surrogates, the Gaussian density is a suitable choice for $\theta(x)$. With this selection, the first $2n$ power moments of $\hat{\rho}(x)$ exist and are finite. Now, we must determine the mean and variance of the Gaussian distribution.

For the DPPM proposed in \cite{wu2023non}, the power moments $\sigma_{1}$ and $\sigma_{2}$ of the reference density $\theta(x)$ can be calculated using (\ref{MomentUpdate}). By choosing $m = \sigma_{1}$ and $\sigma^{2} > \sigma_{2}$ and specifying the density $\theta(x) = \mathcal{N}(m, \sigma^{2})$, we consistently achieve good estimations. The reason behind this is that a relatively large variance $\sigma^{2}$ helps adjust the estimate to densities with multiple peaks (modes).

For the DPBM, we can directly choose the reference density $\theta(x)$ as $\theta(x) = \mathcal{N}(\sigma_{1}, \sigma_{2})$. Due to the additional information provided by the generalized logarithmic moments, we no longer need to choose a relatively larger variance for the prior density.

We first simulate a mixture of Gaussians with two modes,
\begin{equation}
    \rho(x) = \frac{0.5}{\sqrt{2\pi}}e^{-\frac{(x-2)^{2}}{2}} + \frac{0.5}{\sqrt{2\pi}}e^{-\frac{(x+2)^{2}}{2}}.  
\end{equation}
where we select $\theta(x)$ as $\mathcal{N}(0, 5^{2})$. The degree of the polynomial $Q(x)$ is $4$ for both $\hat{\rho}_{m}$ and $\hat{\rho}_{l}$, where  $\hat{\rho}_{m}$ corresponds to DPPM and $\hat{\rho}_{l}$ to DPBM. The highest order of $P(x)$ in $\hat{\rho}_{l}$ is $4$. By Algorithm 2 in \cite{wu2023non}, a non-Gaussian Bayesian filter using only power moments, we obtain $\hat{\rho}_{m} = \theta(x) / Q_{m}(x)$, where $Q_{m}(x) = 4.13 \cdot 10 ^{-2}x^{4} + 5.40 \cdot 10 ^{-5}x^{3} - 4.44 \cdot 10 ^{-1}x^{2} - 3.07 \cdot 10 ^{-4}x + 1.40$. Also we obtain $\hat{\rho_{l}} = \theta(x) \cdot P_{l}(x) / Q_{l}(x)$ by Algorithm \ref{Algo2}, where $Q_{l}(x) = 3.40 \cdot 10 ^{-2}x^{4} - 4.39 \cdot 10 ^{-24}x^{3} - 2.40 \cdot 10 ^{-1}x^{2} - 2.88 \cdot 10 ^{-22}x + 1$ and $P_{l}(x) = 3.19 \cdot 10 ^{-3}x^{4} - 8.39 \cdot 10 ^{-24}x^{3} + 1.60 \cdot 10 ^{-1}x^{2} + 4.49 \cdot 10 ^{-22}x + 1.06$. 

\begin{figure}[htbp]
\centering
\includegraphics[scale=0.35]{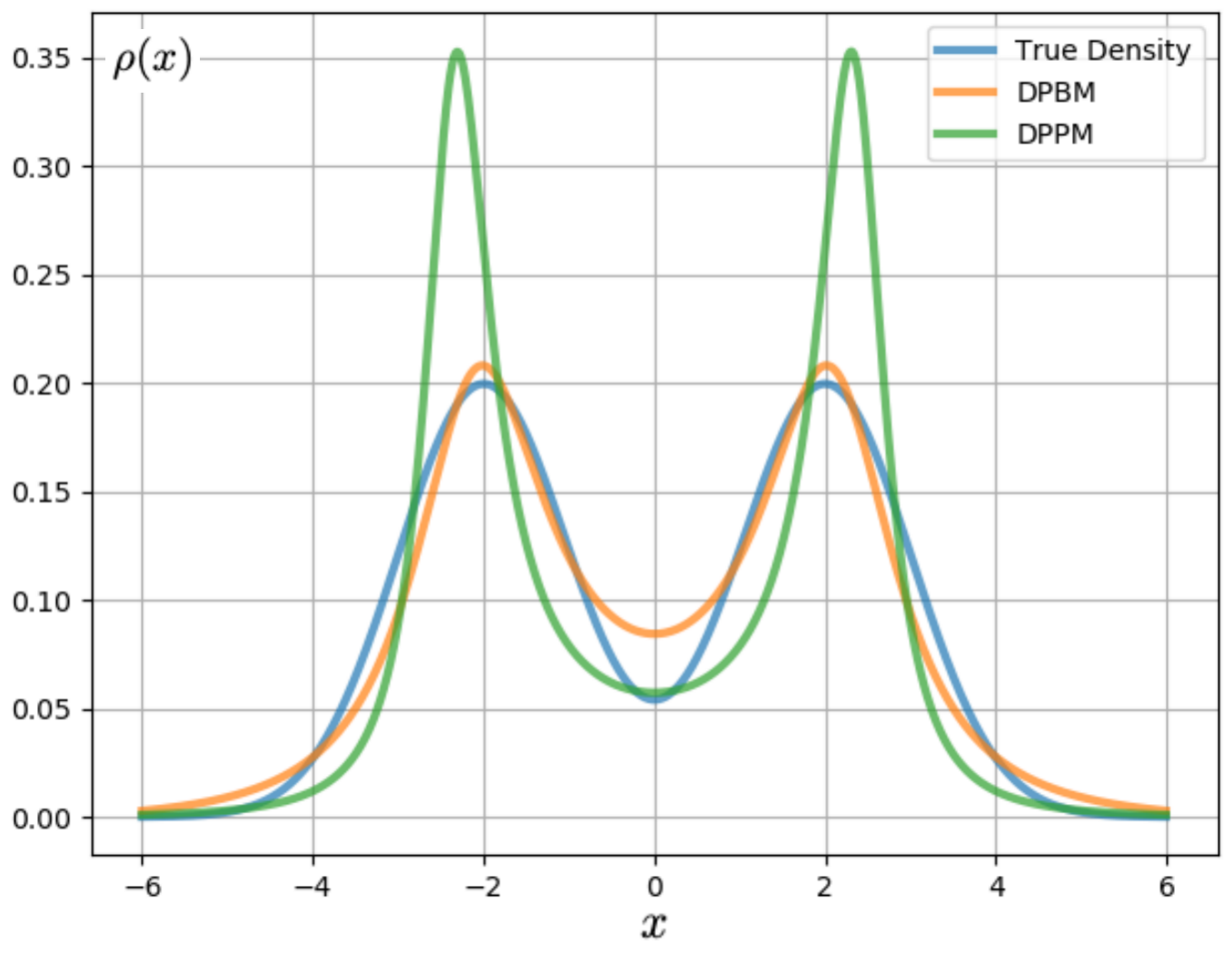}
\centering
\caption{Simulation results of Example 1. The blue curve represents the true prior density function. The green one represents the density estimate using only power moments. And the orange one represents the density estimate using both the power moments and the generalized logarithmic moments.}
\label{fig1}
\end{figure}

The density estimates $\hat{\rho}_{m}$ and $\hat{\rho}_{l}$, along with the true density $\rho$, are depicted in Figure \ref{fig1}. The simulation results clearly demonstrate that incorporating the generalized logarithmic moments into the density surrogate enhances the accuracy of prior density estimation for non-Gaussian Bayesian filtering.

\begin{figure}[htbp]
\centering
\includegraphics[scale=0.35]{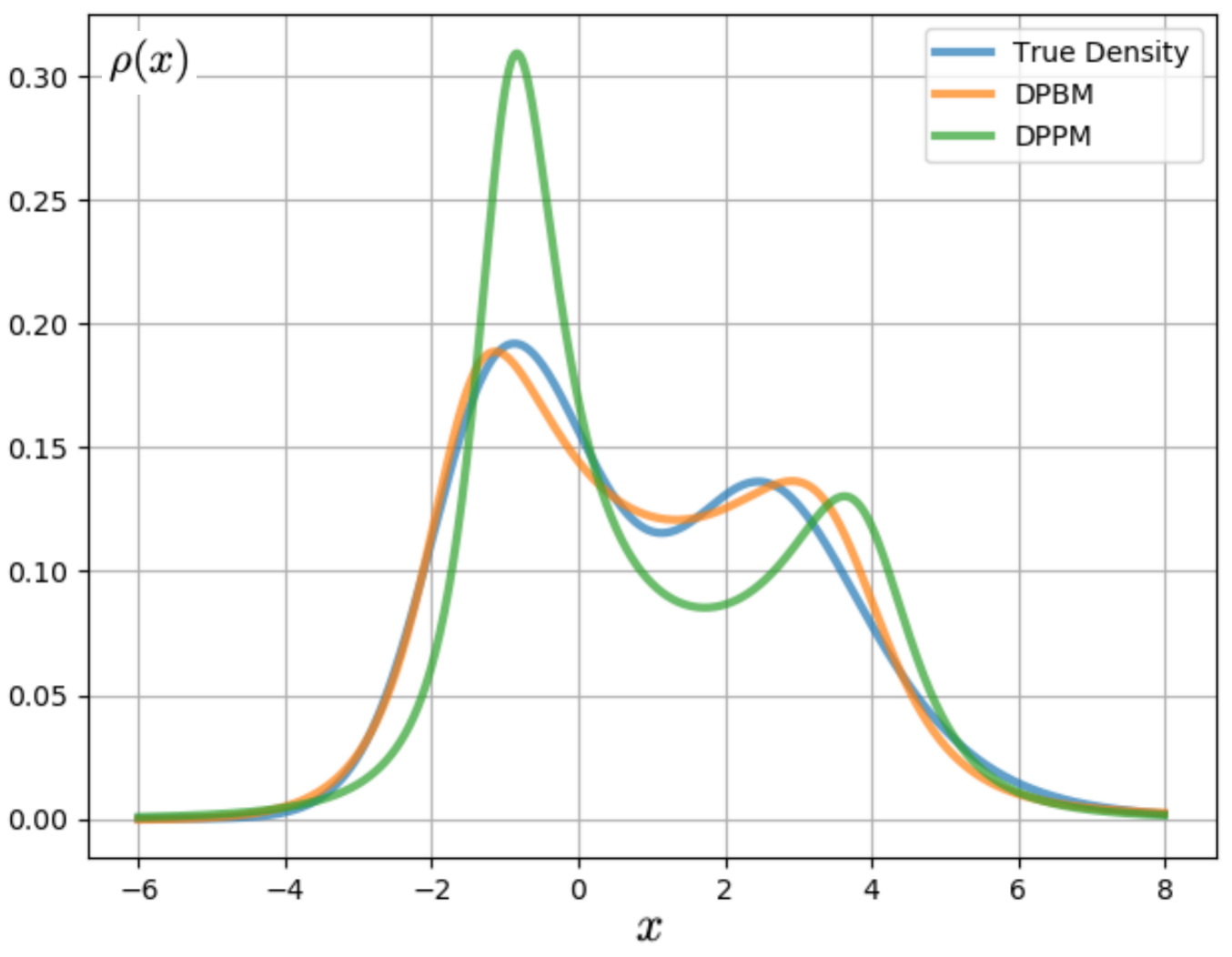}
\centering
\caption{Simulation results of Example 2.}
\label{fig2}
\end{figure}

In the next example, we simulate a mixture of generalized logistic densities, which is known to be challenging to estimate accurately. Specifically, Example 2 represents a mixture of two type-\uppercase\expandafter{\romannumeral1} generalized logistic densities with the probability density function given by
$$
     \rho(x) = \frac{0.4 \cdot 2e^{-x + 2}}{(1 + e^{-x + 2})^{3}} + \frac{0.6 \cdot 3e^{-x - 2}}{(1 + e^{-x - 2})^{4}}.
$$
We choose $\theta(x)$ as $\mathcal{N}(0.90, 5.86^{2})$ as the reference density. For both $\hat{\rho}_{m}$ and $\hat{\rho}_{l}$, we use a degree-4 polynomial $Q(x)$. In $\hat{\rho}_{l}$, the highest order of $P(x)$ is 4. By employing the density surrogates, we obtain $\hat{\rho}_{m} = \frac{\theta(x)}{Q_{m}(x)}$, where $Q_{m}(x) =1.65 \cdot 10^{-2}x^{4} - 9.95 \cdot 10^{-2}x^{3} + 5.27 \cdot 10^{-2}x^{2} + 3.48 \cdot 10^{-1}x + 4 \cdot 10^{-1}$, and $\hat{\rho}_{l} = \frac{\theta(x) \cdot P_{l}(x)}{Q_{l}(x)}$, where $Q_{l}(x) = 1.68 \cdot 10^{-2}x^{4} - 6.82 \cdot 10^{-2}x^{3} - 6.75 \cdot 10^{-2}x^{2} + 3.34 \cdot 10^{-1}x + 1$ and $P_{l}(x) = 7.10 \cdot 10^{-4}x^{4} + 1.75 \cdot 10^{-3}x^{3} - 6.65 \cdot 10^{-2}x^{2} + 9.76 \cdot 10^{-2}x + 2.14$. The simulation results are presented in Figure \ref{fig2}.

It is observed that the estimate $\hat{\rho}_{m}$ is significantly biased compared to the true density. However, by utilizing the density surrogate that incorporates both moments, we achieve an estimate with a considerably reduced error. This outcome is remarkable since the density surrogate only requires 10 parameters and does not rely on any knowledge of $\rho_{x_{t+1} \mid \mathcal{Y}_{t}}(x)$, such as the number of modes or the specific functional form.

In the final example, Example 3, we consider a mixture of two Laplacian densities. The probability density function is defined as follows:
$$
    \rho(x) = 0.3 \cdot e^{-\left| \frac{x - 1}{2}\right|} + 0.7 \cdot e^{-\left| \frac{x + 1}{2}\right|}.  
$$

\begin{figure}[htbp]
\centering
\includegraphics[scale=0.35]{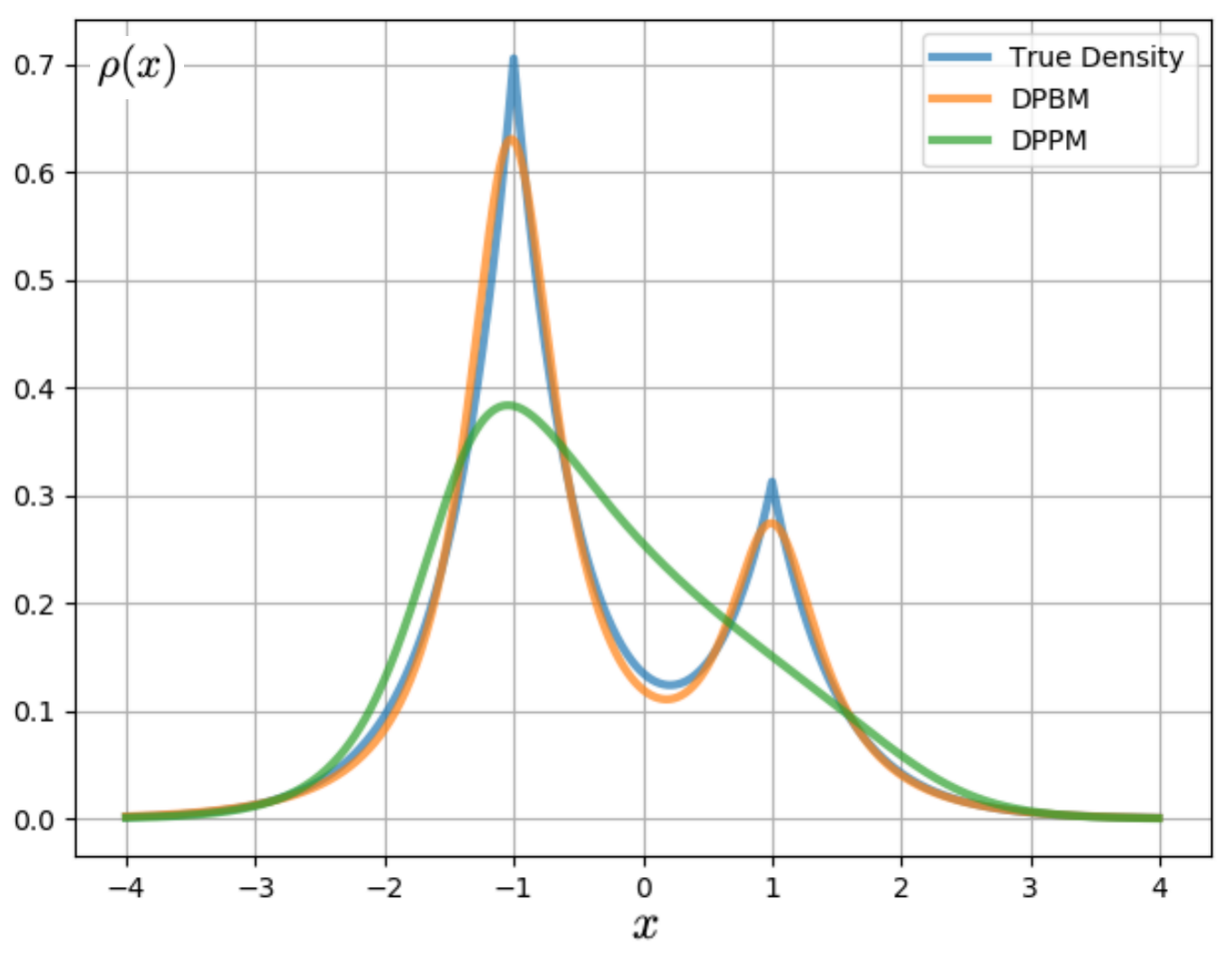}
\centering
\caption{Simulation results of Example 3.}
\label{fig3}
\end{figure}

We select $\theta(x)$ as $\mathcal{N}(-0.4, 1.5^{2})$ for the reference density. The polynomial $Q(x)$ has a maximum order of 4 for both $\hat{\rho}_{m}$ and $\hat{\rho}_{l}$. In $\hat{\rho}_{l}$, the highest order of $P(x)$ is also 4. Utilizing the density surrogates, we obtain $\hat{\rho}_{m} = \frac{\theta(x)}{Q_{m}(x)}$, where $Q_{m}(x) =5.52 \cdot 10^{-2}x^{4} - 7.54 \cdot 10^{-2}x^{3} - 1.69 \cdot 10^{-1}x^{2} + 3.25 \cdot 10^{-1}x + 1.01$, and $\hat{\rho}_{l} = \frac{\theta(x) \cdot P_{l}(x)}{Q_{l}(x)}$, where $Q_{l}(x) = 6.78 \cdot 10^{-1}x^{4} + 5.48 \cdot 10^{-2}x^{3} - 1.11 \cdot x^{2} + 2.39 \cdot 10^{-2}x + 1$ and $P_{l}(x) = 7.00 \cdot 10^{-2}x^{4} + 1.03 \cdot 10^{-1}x^{3} + 6.20 \cdot 10^{-1}x^{2} - 3.13 \cdot 10^{-1}x + 4.67 \cdot 10^{-1}$. The simulation results are depicted in Figure \ref{fig3}.

In this example, it is important to highlight that the density being estimated is not smooth, and it exhibits two distinct sharp modes (peaks) that are in close proximity to each other. When utilizing only the power moments, we observe that the density estimate is unable to accurately capture the two modes, resulting in a poor approximation where only a single peak is represented. However, by incorporating the generalized logarithmic moments into the density surrogate, we significantly enhance the performance of the estimate. The resulting density approximation now successfully captures the presence of the two sharp modes and provides a much-improved representation. 

In the previous numerical examples, we proved that by using the generalized logarithmic moments together with the power moments, the performance of distribution approximation is clearly improved in some situations, compared to merely using the power moments. In the following part of this section, we will give a more illustrative example, where our proposed algorithm is applied to an engineering problem and the performance is compared to several prevailing methods.

We address a robot localization challenge wherein a sensor is designated to measure the distances between the robot and the predefined landmark. We assume that the robot moves along the real line $\mathbb{R}$ with coordinate $x$. The robot's position with respect to coordinate $x$ at time step $t$ is denoted as $x_{t}$. The position of the landmark is denoted as $\check{x}$. In this localization task, the robot undergoes incremental movement, advancing one unit along the positive $x$ direction. Imperfect controls result in deviations from the commanded movement, necessitating consideration of noise in the particle's movements to capture the actual robot movement. The moving distance, with a true value of $1$, is corrupted with an additive Gaussian noise $\mathcal{N}(0, 0.03^{2})$. Furthermore, the distance observation of each landmark is subject to additive noise.

The system and observation equations are
$$
x_{t+1} = x_{t} + 1 + w_{t}
$$
and
$$
z_{t} = x_{t}-\check{x} + v_{t}
$$
respectively. The sign of the signed distance $z_{t}$ is negative when the robot is positioned to the left of the landmark and positive when it is to the right of the landmark.

We assume that the variable $w_{t}$ follows a Gaussian distribution, specifically $\mathcal{N}(0, 0.03^{2})$, accounting for the error in controlling the robot. In prior results, the noises $v_{t}$ was assumed to be Gaussian to obtain a closed form of solution. In this experiment, however, we propose employing the right-skewed Gumbel distribution to validate the proposed algorithm in treating non-Gaussian filtering tasks.

\begin{figure}[htbp]
\centering
\includegraphics[scale=0.35]{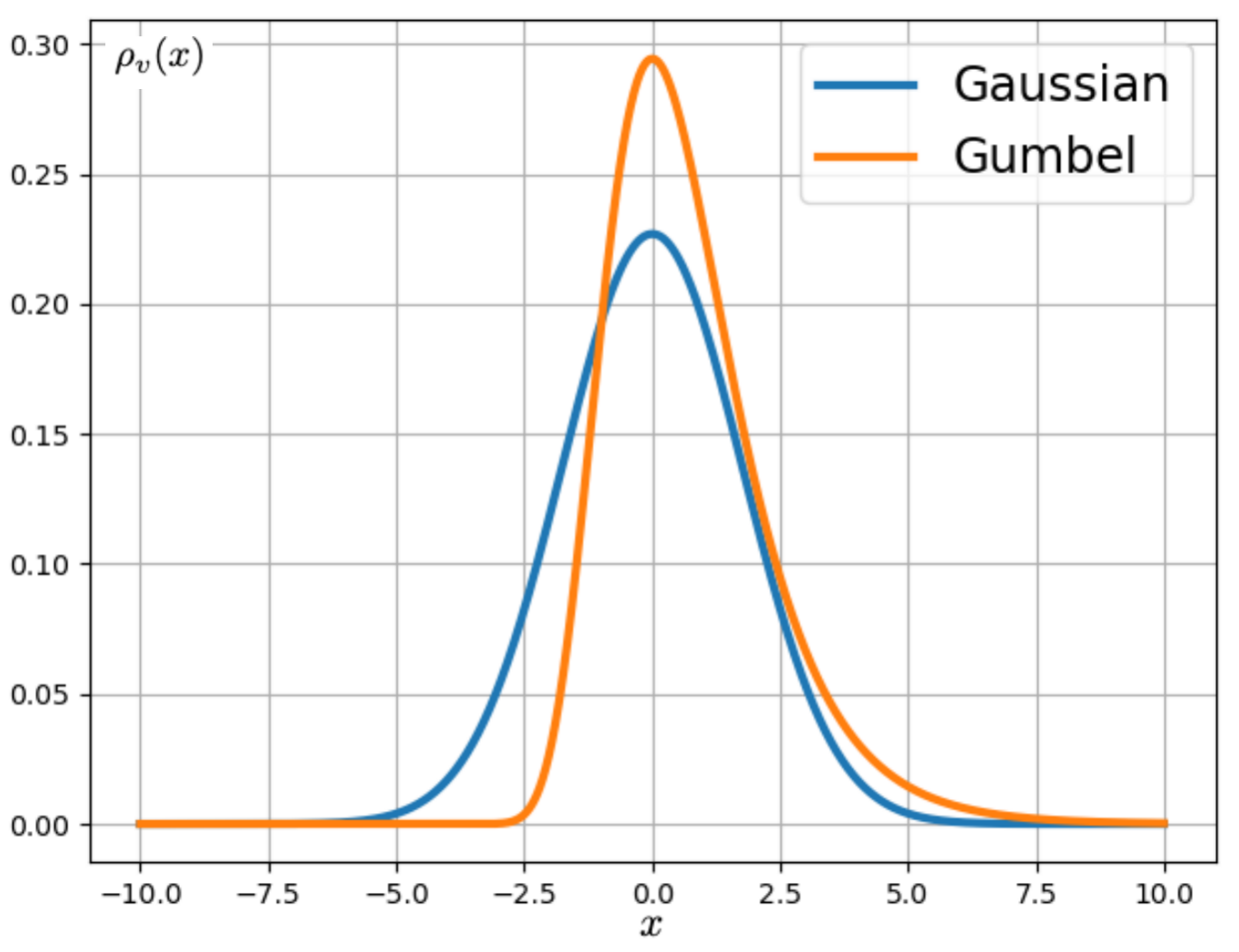}
\centering
\caption{Probability density functions of the Gaussian and the Gumbel distributions.}
\label{fig4}
\end{figure}

The probability density functions of the Gumbel and Gaussian distributions is illustrated in Figure \ref{fig4}. The probability density function of the Gumbel distribution is given by

\begin{equation}
\rho_{v}(x) = 4 e^{-4x-e^{-4x}}.
\label{pv1}
\end{equation}

Meanwhile, the Gaussian distribution, with an identical mean and variance as the Gumbel distribution, has the following probability density function

\begin{equation}
\rho_{v}(x) = \frac{1}{\sqrt{2\pi} \cdot 0.35}e^{\frac{x^{2}}{2 \cdot 0.35^{2}}}.
\label{pv2}
\end{equation}

In the context of this localization task, the use of the asymmetric Gumbel distribution as the model for observation noise poses a significant challenge in selecting an appropriate stochastic filter. Except for DPBM proposed in this paper and the DPPM proposed in \cite{wu2023non}, the particle filters (PF) is the sole feasible option in prevailing methods for carrying out this task, due to the Gumbel distribution. In our simulations, we adopt a sampling-importance resampling (SIR) filter, as described in \cite{chen2003bayesian}. Given that the system equation and the observation equation are both linear, we also adopt the Kalman filter (KF) for this task. However, the KF faces difficulties in handling the Gumbel observation noise. Consequently, we resort to using the Gaussian distribution in \eqref{pv2} as a substitute for \eqref{pv1} to represent the observation noise during the filtering process.

The initialization procedure for the three filters is as follows: The robot initiates its movement from the position $x = -7$, and a landmark is situated at $x = 0$. The distributions of the initial states $x(0)$ for both the DPBM and KF are set to be the Gaussian distribution $\mathcal{N}\left(m_{0}, 1 \right)$, where the mean $m_{0}$ is drawn from the Gaussian distribution $\mathcal{N}\left(-7, 1 \right)$. The states of the $5000$ particles in the PF are i.i.d samples drawn from the uniform distribution $U \left( [-8, 8] \right)$, which aims to cover a broader range of potential locations. The additive noise in the distance observation follows the Gumbel distribution in \eqref{pv1}. The DPBM utilizes power moments and generalized logarithmic moments both up to the fourth order to estimate the density surrogates.

\begin{figure}[htbp]
\centering
\includegraphics[scale=0.35]{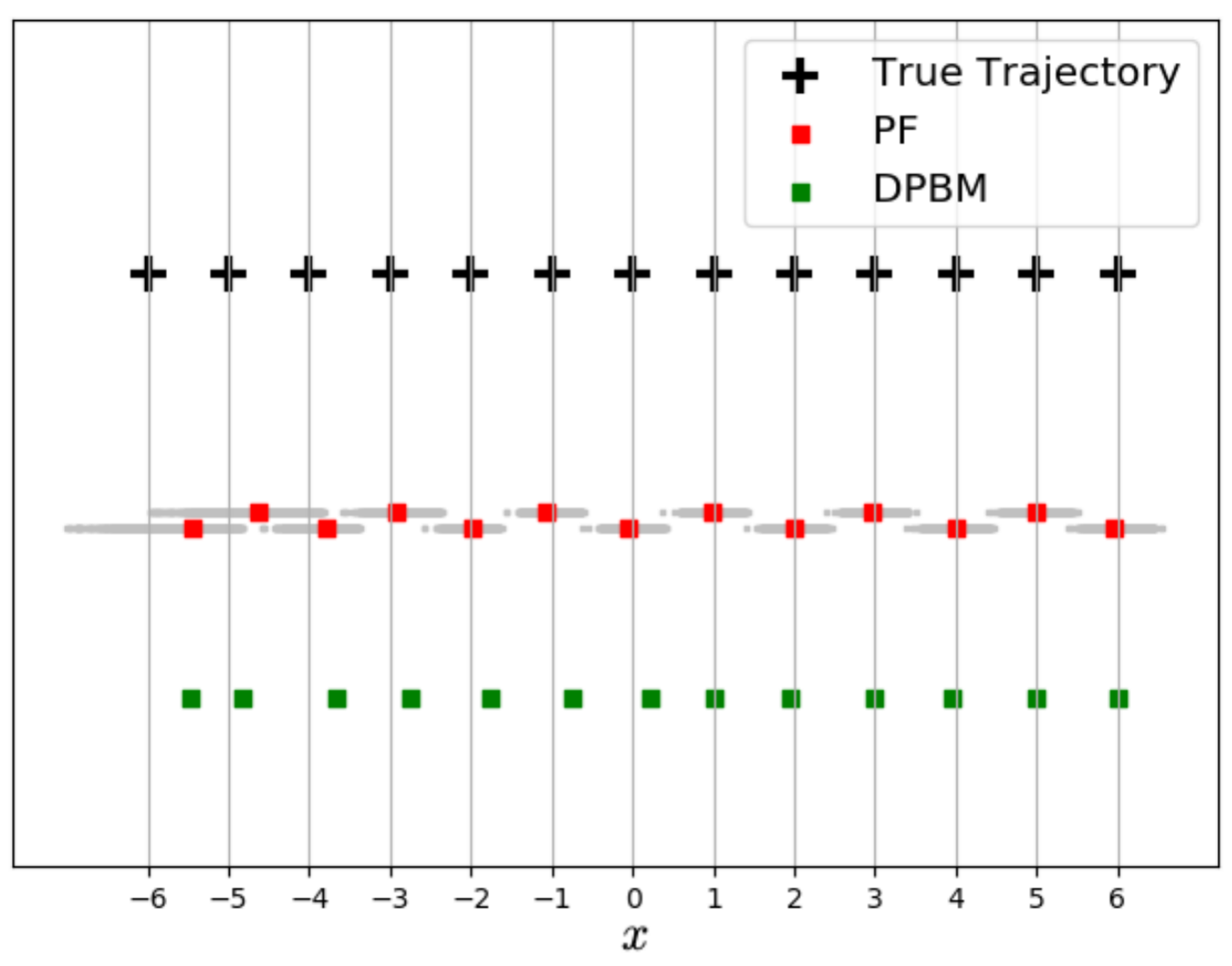}
\centering
\caption{A sample localization process.}
\label{fig5}
\end{figure}

\begin{figure}[htbp]
\centering
\includegraphics[scale=0.35]{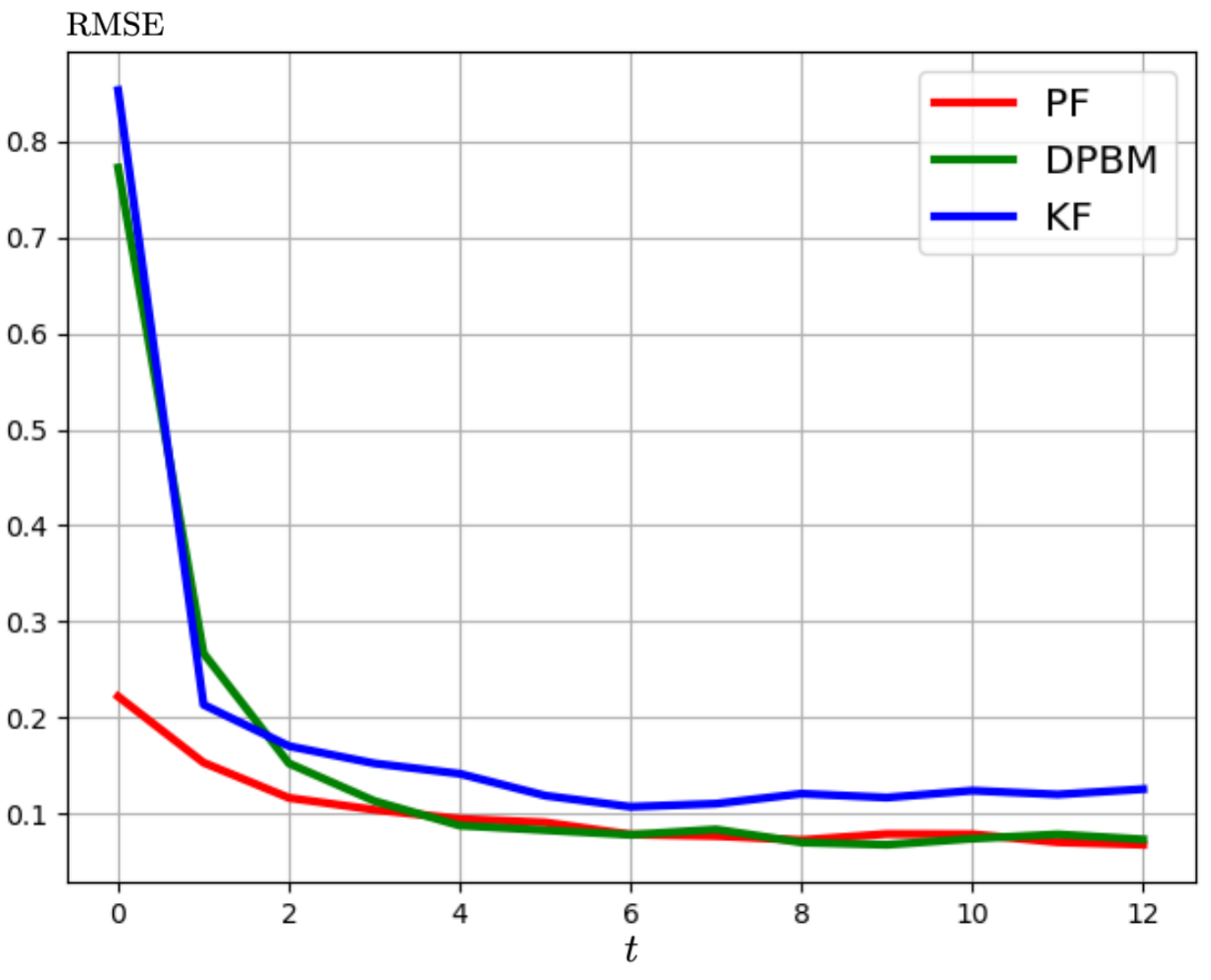}
\centering
\caption{RMSE as a function of time step $t = 1, \cdots, 13$ of $50$ Monte-Carlo simulations for KF, PF and DPBM.}
\label{fig6}
\end{figure}

Figure \ref{fig5} illustrates a sample robot localization process along the x-axis, showcasing estimation results by PF and DPBM. The black crosses represent the true trajectory of the robot. The red and green dots represent the location estimates by the particle filter and our proposed DPBM filter. The gray dots represent the particles of the particle filter at each time step. We emphasize that the robot moves along the x-axis, even though the estimates and the particles of PF at different time steps are drawn in two levels, which aims to show the locations of the particles better. Notably, the location estimates by the DPBM converge to the true locations, while the PF particle states also converge to the correct positions. Figure~\ref{fig6} presents the root mean square error (RMSE) curves for 50 Monte-Carlo simulations of DPBM, PF, and KF. As the state estimates converge, the RMSE of DPBM is almost equal to that of PF. We note that the RMSE of KF is significantly larger due to the use of a Gaussian distribution as an alternative for the true Gumbel distribution. The biased observation noise model causes obvious performance degradation of filtering.

From an RMSE standpoint, the DPBM does not outperform the PF, but a notable drawback of the Particle filter is its requirement to store massive amount of data. For instance, in this simulation, the state of each particle consists of two parameters, namely its position and weight, resulting in a need for 10,000 parameters to characterize the system state density. In contrast, the DPBM only requires $9$ parameters for this task {($4$ parameters for the numerator of \eqref{MiniForm} and $5$ parameters for the denominator of \eqref{MiniForm})}, offering a more compact representation of the density function.

Additionally, considering the time consumption, the PF outperforms the DPBM in the execution time. In this example, each filtering iteration takes an average execution time of 1.08 seconds on a 2.5 GHz Intel Core i7 CPU. While this may be relatively long compared to the PF execution time, it remains manageable for applications with less sensitivity to the execution time. {For example, intraday, daily, weekly, monthly and annual data are widely used in analyzing the futures market \cite{quintana2010futures}, where filters serve as approaches for this task \cite{manoliu2002energy, fileccia2018particle}. Our proposed filtering scheme is very suitable for this task.} Moreover, the optimization in each filtering step is convex with the solution proved to exist and be unique. It makes the execution for each filtering step to be predictable, which is a clear advantage of DPBM.

The comparison presented above solely examines the errors in the first-order statistics of the filtering results, potentially obscuring a distinct advantage of the proposed scheme over PF. With the proposed filtering approach, the prior and posterior probability density functions for each filtering iteration are continuous rational functions. Consequently, we can easily obtain the probability value for any arbitrary $x \in \mathbb{R}$. Conversely, this is not feasible with PF due to the discrete nature of particles. 

Theoretically speaking, the proposed filter is a consistent one, which is the best filter for the cases where exact filters don't exist. The simulation results further confirm the validity of the proposed continuous approximation, which utilizes power and generalized logarithmic moments, for the intractable integral. The proposed density parametrization and the particle-representation of the particle filter represent two types of methods to treat the intractable integral, which provide consistent continuous and discrete estimators of the integral respectively.

\section{Conclusions}

A Bayesian filter based on density parametrization using both the power moments and the generalized logarithmic moments of the densities is developed in this paper. We propose a convex optimization scheme to uniquely determine a rational density with exactly the specified power and generalized logarithmic moments, rather than estimating the parameters of a prespecified density model (such as Gaussian or Student's t) by minimizing the difference characterized by a norm, like the traditional method of moments. The map from the parameters of the proposed density surrogate to the power and generalized logarithmic moments is proved to be diffeomorphic, which reveals the fact that the parameters can be uniquely determined by the two types of moments. Furthermore, we provide the statistical property, together with numerical simulations to validate the proposed density estimator. By the results of the numerical simulations, we observe that the performance of density estimation using the proposed algorithm is quite satisfactory, which is a clear improvement as compared to that of the density surrogate using only power moments in our previous paper \cite{wu2023non}. The presented filter is employed in a robot localization task, enabling a comprehensive performance comparison with various established filtering schemes. It is noteworthy that the Root Mean Square Error (RMSE) exhibited by the proposed filter is comparable to that of the commonly used particle filter. Despite the fact that the execution time for each filtering step in the proposed algorithm is comparatively longer than that of the particle filter, the proposed algorithm effectively mitigates the necessity for the particle filter to store states of massive particles.

\begin{ack}                               
We would like to express our sincere appreciation to the anonymous reviewers and the associate editor for their valuable time, expertise, and constructive feedback provided during the review process of this journal paper. Their insightful comments and suggestions have significantly contributed to improving the quality and clarity of our work.
\end{ack}

\bibliographystyle{elsarticle-num}
\bibliography{autosam}

\begin{appendix}
\section{Proof of Theorem \ref{MomentError}}
\label{Appendix_asymptotic}
\begin{proof} 
It is proved in \cite{wu2023non} that, for a sufficiently large $n$,
\begin{equation}
\begin{aligned}
    & \mathbb{E}\left ( x_{t}^{k}|\mathcal{Y}_{t} \right ) \approx \mathbb{E}\left ( \hat{x}_{t}^{k}|\mathcal{Y}_{t} \right )\\
    \text{and} \quad & \mathbb{E}\left ( x_{t+1}^{k}|\mathcal{Y}_{t} \right ) \approx \mathbb{E}\left ( \hat{x}_{t+1}^{k}|\mathcal{Y}_{t} \right )
\label{MomentApprox}
\end{aligned}
\end{equation}
for $k = 1, \cdots, 2n$,  for $\rho_{x_t|\mathcal{Y}_t}, \rho_{x_{t+1}|\mathcal{Y}_t} \in \mathcal{SG}$, where $\mathcal{SG}$ denotes the space of all sub-Gaussian distributions. Sub-Gaussian distributions are those whose tails are dominated by the tails of a Gaussian distribution, i.e., decay at least as fast as a Gaussian. Moreover,
\begin{equation}
\begin{aligned}
    & \lim_{n \rightarrow +\infty}\mathbb{E}\left ( \hat{x}_{t}^{k}|\mathcal{Y}_{t} \right ) = \mathbb{E}\left ( x_{t}^{k}|\mathcal{Y}_{t} \right )\\
    \text{and} \quad & \lim_{n \rightarrow +\infty}\mathbb{E}\left ( \hat{x}_{t+1}^{k}|\mathcal{Y}_{t} \right ) = \mathbb{E}\left ( x_{t+1}^{k}|\mathcal{Y}_{t} \right )
\label{MomentInf}
\end{aligned}
\end{equation}
for $k = 1, \cdots, +\infty$. Now it remains to analyze
$$
\begin{aligned}
    & \mathbb{E}^{\log}\left ( x_{t}^{k}|\mathcal{Y}_{t} \right ) - \mathbb{E}^{\log}\left ( \hat{x}_{t}^{k}|\mathcal{Y}_{t} \right )\\
    \text{and} \quad & \mathbb{E}^{\log}\left ( x_{t+1}^{k}|\mathcal{Y}_{t} \right ) - \mathbb{E}^{\log}\left ( \hat{x}_{t+1}^{k}|\mathcal{Y}_{t} \right )
\end{aligned}
$$

for $k = 1, \cdots, 2n$. We note that $\mathbb{E}^{\log}\left ( x_{1}^{k}|\mathcal{Y}_{0} \right ) = \mathbb{E}^{\log}\left ( \hat{x}_{1}^{k}|\mathcal{Y}_{0} \right )$ after the first time update, i.e.,
\begin{equation}
\begin{aligned}
    \int_{\mathbb{R}} x^{k} \theta(x) \left[ \log\left ( \rho_{x_{1} \mid \mathcal{Y}_{0}}\right) - \log\left(\hat \rho_{x_{1} \mid \mathcal{Y}_{0}} \right )\right]dx = 0.
\label{x1y0}
\end{aligned}
\end{equation}
Meanwhile, we can write the generalized logarithmic moment terms of $\rho_{x_{1} \mid \mathcal{Y}_{1}}$ as
\begin{equation*}
    \mathbb{E}^{\log}\left ( x_{1}^{k}|\mathcal{Y}_{1} \right ) = \int_{\mathbb{R}} x^{k} \theta(x)\log\left[ \rho_{\epsilon_{1}}\left(y_{1}-h_{1} x\right) \rho_{x_{1} \mid \mathcal{Y}_{0}}(x)\right]dx
\end{equation*}
for $k = 1, \cdots, 2n$, and those of $\hat \rho_{x_{1} \mid \mathcal{Y}_{1}}$ as
\begin{equation*}
    \mathbb{E}^{\log}\left ( \hat{x}_{1}^{k}|\mathcal{Y}_{1} \right ) = \int_{\mathbb{R}} x^{k} \theta(x)\log\left[\rho_{\epsilon_{1}}\left(y_{1}-h_{1} x\right) \hat \rho_{x_{1} \mid \mathcal{Y}_{0}}(x)\right]dx
\end{equation*}
for $k = 1, \cdots, 2n$. Therefore by \eqref{x1y0} we have,
\begin{equation*}
\begin{aligned}
    & \mathbb{E}^{\log}\left ( x_{1}^{k}|\mathcal{Y}_{1} \right ) - \mathbb{E}^{\log}\left ( \hat{x}_{1}^{k}|\mathcal{Y}_{1} \right )\\
    = & \int_{\mathbb{R}} x^{k} \theta(x) \left[ \log\left ( \rho_{x_{1} \mid \mathcal{Y}_{0}}\right) - \log\left(\hat \rho_{x_{1} \mid \mathcal{Y}_{0}} \right )\right]dx\\
    = & 0
\end{aligned}
\end{equation*}
for $k = 1, \cdots, 2n$. Then we have
\begin{equation*}
    \mathbb{E}^{\log}\left ( x_{1}^{k}|\mathcal{Y}_{1} \right ) = \mathbb{E}^{\log}\left ( \hat{x}_{1}^{k}|\mathcal{Y}_{1} \right ),\quad k = 1, \cdots, 2n.
\end{equation*}

Moreover, by \eqref{Prediction}, we have
\begin{equation*}
\begin{aligned}
    & \mathbb{E}^{\log}\left ( x_{2}^{k}|\mathcal{Y}_{1} \right ) - \mathbb{E}^{\log}\left ( \hat{x}_{2}^{k}|\mathcal{Y}_{1} \right )\\
    = & \int_{\mathbb{R}} x^{k} \theta(x) \log \int_{\mathbb{R}} \rho_{x_{1} \mid \mathcal{Y}_{1}}\left(\frac{\varepsilon}{f_{1}}\right) \rho_{\eta_{1}}(x-\varepsilon) d\varepsilon dx\\
    - & \int_{\mathbb{R}} x^{k} \theta(x) \log \int_{\mathbb{R}} \hat{\rho}_{x_{1} \mid \mathcal{Y}_{1}}\left(\frac{\varepsilon}{f_{1}}\right) \rho_{\eta_{1}}(x-\varepsilon) d\varepsilon dx\\
    = & \int_{\mathbb{R}}f_{1} x^{k} \theta(x) \log \int_{\mathbb{R}} \rho_{x_{1} \mid \mathcal{Y}_{1}}\left(\omega\right) \rho_{\eta_{1}}(x-f_{1}\omega) d\omega dx\\
    - & \int_{\mathbb{R}}f_{1} x^{k} \theta(x) \log \int_{\mathbb{R}} \hat{\rho}_{x_{1} \mid \mathcal{Y}_{1}}\left(\omega\right) \rho_{\eta_{1}}(x-f_{1}\omega) d\omega dx
\end{aligned}
\end{equation*}
for $k = 1, \cdots, 2n$. 
We note that $\rho_{\eta_{1}}\left(x - f_{1}\omega \right)$ is analytic almost everywhere. Assume $\rho_{\eta_{1}}\left(x - f_{1}\omega \right)$ is analytic at point $x_{0}$, then it is feasible for us to write the Taylor series at this point. Without loss of generality, we take $x_{0} = f_{1}\omega$, then we have
$$
\begin{aligned}
    & \rho_{\eta_{1}}\left(x - f_{1}\omega \right)\\
    = & \sum_{i = 0}^{+\infty} \frac{\rho^{(i)}_{\eta_{1}}\left(0\right)}{i!}\left(x - f_{1}\omega\right)^{i}\\
    = & \sum_{i = 0}^{+\infty} \sum_{j = 0}^{i} \left ( \begin{matrix} i\\j\end{matrix} \right )\frac{(-f_{1})^{j}\rho^{(i)}_{\eta_{1}}\left(0\right)}{i!}\omega^{j}x^{i-j}\\
\end{aligned}
$$

Since all power moments and generalized logarithmic moments of $x_{1}$ and $\hat{x}_{1}$ exist and are finite, we have \eqref{longeq1}. 

\begin{figure*}[t]
\begin{equation}
\begin{aligned}
    & \mathbb{E}^{\log}\left ( x_{2}^{k}|\mathcal{Y}_{1} \right ) - \mathbb{E}^{\log}\left ( \hat{x}_{2}^{k}|\mathcal{Y}_{1} \right )\\
    = & \int_{\mathbb{R}} x^{k} \theta(x) \log \int_{\mathbb{R}} \rho_{x_{1} \mid \mathcal{Y}_{1}}\left(\omega \right) \sum_{i = 0}^{+\infty} \sum_{j = 0}^{i}\left ( \begin{matrix} i\\j\end{matrix} \right ) \frac{(-f_{1})^{j}\rho^{(i)}_{\eta_{1}}\left(0\right)}{i!}\omega^{j}x^{i-j} d\omega dx\\
    - & \int_{\mathbb{R}} x^{k} \theta(x) \log \int_{\mathbb{R}} \hat{\rho}_{x_{1} \mid \mathcal{Y}_{1}}\left(\omega\right) \sum_{i = 0}^{+\infty} \sum_{j = 0}^{i}\left ( \begin{matrix} i\\j\end{matrix} \right ) \frac{(-f_{1})^{j}\rho^{(i)}_{\eta_{1}}\left(0\right)}{i!}\omega^{j}x^{i-j} d\omega dx\\
    = & \int_{\mathbb{R}} x^{k} \theta(x) \log \left( \sum_{i = 0}^{+\infty} \sum_{j = 0}^{i}\left ( \begin{matrix} i\\j\end{matrix} \right )x^{i-j} \int_{\mathbb{R}} \rho_{x_{1} \mid \mathcal{Y}_{1}}\left(\omega\right) \frac{(-f_{1})^{j}\rho^{(i)}_{\eta_{1}}\left(0\right)}{i!}\omega^{j} d\omega \right) dx\\
    - & \int_{\mathbb{R}} x^{k} \theta(x) \log \left( \sum_{i = 0}^{+\infty} \sum_{j = 0}^{i}\left ( \begin{matrix} i\\j\end{matrix} \right )x^{i-j} \int_{\mathbb{R}} \hat{\rho}_{x_{1} \mid \mathcal{Y}_{1}}\left(\omega\right) \frac{(-f_{1})^{j}\rho^{(i)}_{\eta_{1}}\left(0\right)}{i!}\omega^{j} d\omega \right) dx\\
    = & \int_{\mathbb{R}} x^{k} \theta(x) \log \left( \sum_{i = 0}^{+\infty} \sum_{j = 0}^{i}\left ( \begin{matrix} i\\j\end{matrix} \right )x^{i-j} \frac{(-f_{1})^{j}\rho^{(i)}_{\eta_{1}}\left(0\right)}{i!}\mathbb{E}\left ( x_{1}^{j}|\mathcal{Y}_{1} \right ) \right) dx\\
    - & \int_{\mathbb{R}} x^{k} \theta(x) \log \left( \sum_{i = 0}^{+\infty} \sum_{j = 0}^{i}\left ( \begin{matrix} i\\j\end{matrix} \right )x^{i-j} \frac{(-f_{1})^{j}\rho^{(i)}_{\eta_{1}}\left(0\right)}{i!}\mathbb{E}\left ( \hat{x}_{1}^{j}|\mathcal{Y}_{1} \right ) \right) dx\\
\end{aligned}
\label{longeq1}
\end{equation}
\hrulefill
\vspace*{4pt}
\end{figure*}

By \eqref{longeq1}, we note that $\mathbb{E}^{\log}\left ( x_{2}^{k}|\mathcal{Y}_{1} \right ) - \mathbb{E}^{\log}\left ( \hat{x}_{2}^{k}|\mathcal{Y}_{1} \right )$ tends to zero as $n\to \infty$ by \eqref{MomentApprox}. By properly selecting a sufficient large $n$, we have
\begin{equation*}
    \mathbb{E}^{\log}\left ( x_{2}^{k}|\mathcal{Y}_{1} \right ) \approx \mathbb{E}^{\log}\left ( \hat{x}_{2}^{k}|\mathcal{Y}_{1} \right ),\quad k = 1, \cdots, 2n,
\end{equation*}

Similarly we can prove
\begin{equation*}
    \mathbb{E}^{\log}\left ( x_{t}^{k}|\mathcal{Y}_{t} \right ) \approx \mathbb{E}^{\log}\left ( \hat{x}_{t}^{k}|\mathcal{Y}_{t} \right ),\quad k = 1, \cdots, 2n,
\end{equation*}
and
\begin{equation*}
    \mathbb{E}^{\log}\left ( x_{t+1}^{k}|\mathcal{Y}_{t} \right ) \approx \mathbb{E}^{\log}\left ( \hat{x}_{t+1}^{k}|\mathcal{Y}_{t} \right ),\quad k = 1, \cdots, 2n,
\end{equation*}
as claimed.
\end{proof}

\section{Connectivity of $\mathcal{M}_{2n}({\sigma})$}
\label{Append_xi}
It is nontrivial to prove that the set of all feasible $(p_{1}, \cdots, p_{2n})$ is path-connected given $\mathcal{M}_{2n}({\sigma})$. 

From the view of optimization, if the feasible $(p_{1}, \cdots, p_{2n})$ fall into several disjoint sets, it is difficult to achieve the global optimum. In this appendix, we prove the connectivity of $\mathcal{M}_{2n}({\sigma})$.

We first prove that the map sending $(p_{1}, \cdots, p_{2n})$ to $\xi \in \mathcal{C}_{2n}$ is a diffeomorphism.

It is obvious that given a $P$, there exists a unique $\xi$. Now we need to prove that given a generalized logarithmic moment sequence $\xi$, there exists a unique $P$. Here we prove this by contradiction. Assume $\frac{P\left ( x \right )}{Q\left ( x \right )}\theta$ and $\frac{{P'}\left ( x \right )}{Q\left ( x \right )}\theta$ correspond to an identical $\xi$ where $P\left ( x \right ) \neq P'\left ( x \right )$, i.e.,
\begin{equation}
    \int_{\mathbb{R}} \theta x^{i} \log \frac{P}{Q}\theta d x=\int_{\mathbb{R}} \theta x^{i} \log \frac{{P'}}{Q}\theta d x=\xi_{i}
\end{equation}
for $i = 0, 1, \cdots, 2n$ (specifically, $\xi_{0}$ is confined to be zero).

Therefore we have
\begin{equation}
    \int_{\mathbb{R}} x^{i}\theta(x) \log \frac{P(x)\theta(x)}{{P'}(x)\theta(x)} d x = 0, \quad i=0,1, \ldots, 2n.
\end{equation}

As $P(x)\theta(x)$ and ${P'}(x)\theta(x)$ are both positive functions, 
\begin{equation}
\begin{aligned}
    & \sum_{i = 0}^{2n}p_{i}\int_{\mathbb{R}} x^{i}\theta(x) \log \frac{P(x)\theta(x)}{{P'}(x)\theta(x)} d x\\
    & = \int_{\mathbb{R}} \sum_{i = 0}^{2n}p_{i} x^{i}\theta(x) \log \frac{P(x)\theta(x)}{{P'}(x)\theta(x)} d x\\
    & = \mathbb{KL}(P\theta, P'\theta)\\
    & = 0.
\end{aligned}
\label{KLPP}
\end{equation}
However, \eqref{KLPP} does not necessarily lead to the fact that $P(x) = P'(x)$, since $\int_{\mathbb{R}}P(x)\theta(x)dx$ and $\int_{\mathbb{R}}P'(x)\theta(x)dx$ are not necessarily equal to $1$. We now prove the following lemma.

\begin{lemma}
    For any arbitrary $P, P' \in \mathcal{S}_{2n}$, $\mathbb{KL}(P\theta, P'\theta) = 0$ if and only if $P \equiv P'$, i.e.,
$$
\left ( p_{1}, \cdots, p_{2n} \right ) = \left ( p'_{1}, \cdots, p'_{2n} \right ).
$$
\label{LemmaB1}
\end{lemma}

\begin{proof}
The negative Kullback-Leibler distance reads
    \begin{equation}
    \begin{aligned}
        & -\mathbb{KL}(P\theta, P'\theta)\\
        = & \int_{\mathbb{R}} P(x)\theta(x) \log \frac{P'(x)\theta(x)}{P(x)\theta(x)} d x\\
        \leq & \int_{\mathbb{R}} P(x)\theta(x) \left( \frac{P'(x)\theta(x)}{P(x)\theta(x)} - 1 \right)dx\\
        = & \int_{\mathbb{R}} \left(P'(x)\theta(x) - P(x)\theta(x) \right)dx\\
        = & 0,
    \end{aligned}
    \end{equation}
with the equality achieved if and only if $P'(x)\theta(x) \equiv P'(x)\theta(x)$, which completes the proof.
\end{proof}

By Lemma \ref{LemmaB1}, we have that $P \equiv P'$ given a generalized logarithmic moment sequence $\xi$. This contradicts our assumption. We have that the map sending $\left ( p_{1}, \cdots, p_{2n} \right )$ to $( \xi_{1}, \cdots, \xi_{2n})$ is a bijection. And because the map and the inverse map are both differentiable, we have that the map is a diffeomorphism. Therefore, $\mathcal{M}_{2n}(\sigma)$ is diffeomorphic to $\mathcal{M}_{2n}(Q)$, which is again diffeomorphic to $\mathbb{R}^{2n}$ and is then path-connected.  

\section{Connectivity of $\mathcal{M}_{2n}({\xi})$}
\label{Append_sigma}
We will prove that the $2n$-manifold $(q_{1}, \cdots, q_{2n})$ is path-connected given $\mathcal{M}_{2n}({\xi})$. First we prove that the map sending $(q_{1}, \cdots, q_{2n})$ to $\sigma \in \mathcal{R}_{2n}$ is a diffeomorphism.

It is obvious that given a pair of parameters $(q_{1}, \cdots, q_{2n})$, there exists a unique $\sigma$. Now we need to prove that given a specific power moment sequence $\sigma$, there exists a unique $(q_{1}, \cdots, q_{2n})$. Again we prove by contradiction. Assume $\frac{P\left ( x \right )}{Q\left ( x \right )}$ and $\frac{{P}\left ( x \right )}{Q'\left ( x \right )}$ have the identical $\sigma$, i.e.,
\begin{equation}
    \int_{\mathbb{R}} x^{i} \frac{P(x)}{Q(x)}\theta(x) dx=\int_{\mathbb{R}} x^{i} \frac{{P}(x)}{Q'(x)}\theta(x) d x=\sigma_{i},
\end{equation}
for $i=0,1, \ldots, 2n$.

Then we have
\begin{equation}
    \int_{\mathbb{R}} x^{i} \frac{P(x)\left ( Q(x)-Q'(x) \right )}{Q(x)Q'(x)} \theta(x) d x = 0,
\end{equation}
for $i=0,1, \ldots, 2n$.  Defining $(\tilde{q}_i)$ by 
\begin{equation}
    Q(x) - Q'(x) = \sum_{i=0}^{2n}\tilde{q}_{i}x^{i},
\end{equation}
we have
$$
\begin{aligned}
& \sum_{i=0}^{2n}\tilde{q}_{i}\int_{\mathbb{R}} x^{i} \frac{P(x)\left ( Q(x)-Q'(x) \right )}{Q(x)Q'(x)} \theta(x) d x\\
= & \int_{\mathbb{R}} \sum_{i=0}^{2n}\tilde{q}_{i} x^{i} \frac{P(x)\left ( Q(x)-Q'(x) \right )}{Q(x)Q'(x)} \theta(x) d x\\
= & 0,
\end{aligned}
$$
and consequently
\begin{equation}
    \int_{\mathbb{R}}\sum_{i=0}^{2n}P(x) \cdot \frac{\left ( \sum_{i=0}^{2n}\tilde{q}_{i}x^{i} \right )^{2}}{Q(x)Q'(x)}\theta(x)dx = 0 .
\end{equation}

Since $P(x), Q(x)$ and $Q'(x)$ are all positive, we have
\begin{equation}
    \sum_{i=0}^{2n}\tilde{q}_{i}x^{i} \equiv 0,
\end{equation}
i.e.
\begin{equation}
    \tilde{q}_{i} = 0, \quad i = 0, 1, \cdots, 2n
\end{equation}

This contradicts our assumption. Therefore we can conclude that the map sending $\left ( q_{1}, \cdots, q_{2n} \right )$ to $\left ( \sigma_{1}, \cdots, \sigma_{2n} \right )$ is a bijection. And because the map and the inverse map are both differentiable, we have that the map is a diffeomorphism. Since $\mathcal{M}_{2n}(\xi)$ and $\mathcal{M}_{2n}(P)$ are both differentiable, they are diffeomorphic. Then we have that $\mathcal{M}_{2n}(\xi)$ is smooth and path-connected, because $\mathcal{M}_{2n}(Q)$ also has these two properties.
\end{appendix}

\begin{wrapfigure}{l}{20mm} 
\includegraphics[width=1in,height=1.25in,clip,keepaspectratio]{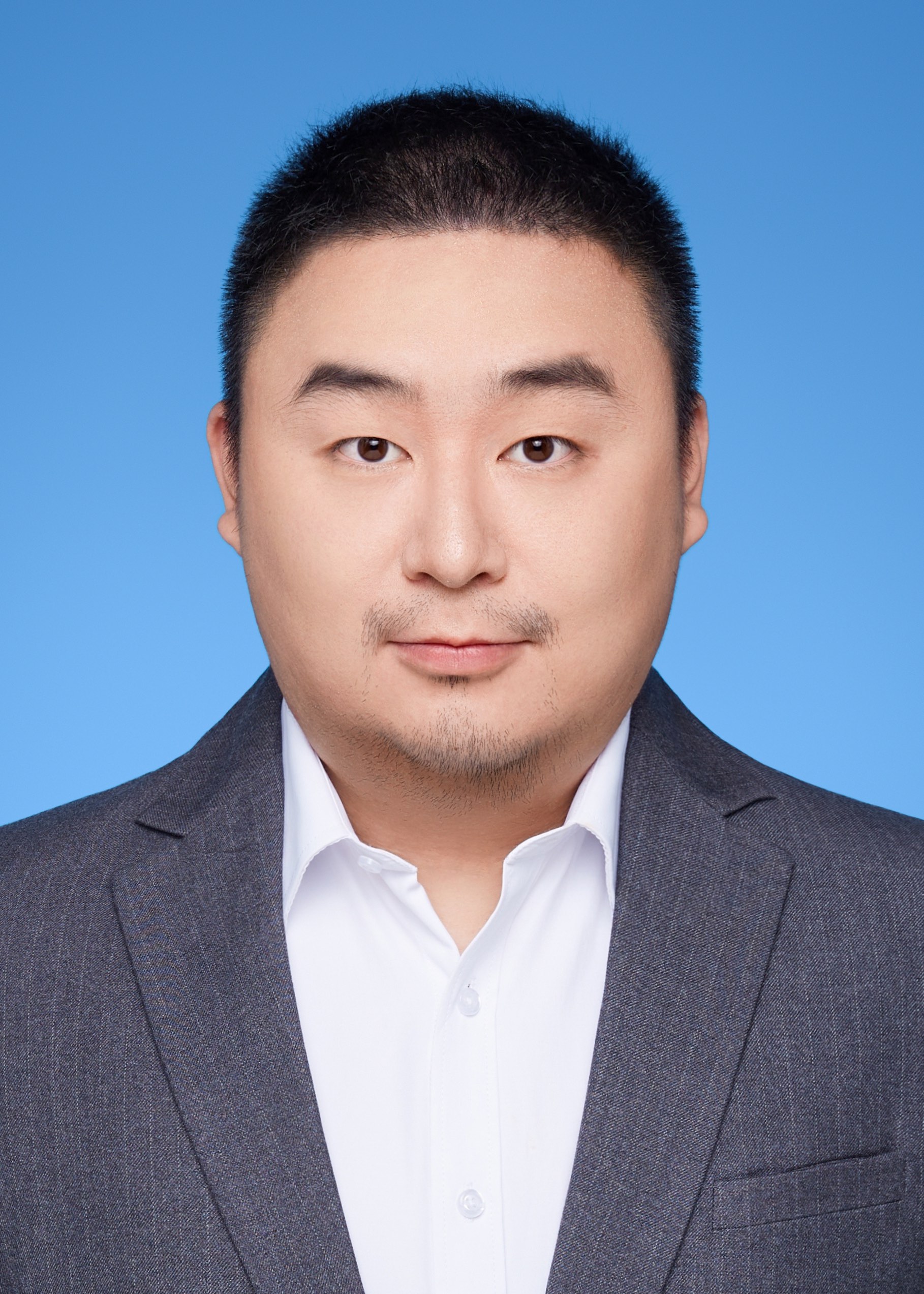}
  \end{wrapfigure}\par
  \textbf{Guangyu Wu} received the B.E. degree from Northwestern Polytechnical University, Xi'an, China, in 2013, and two M.S. degrees, one in control science and engineering from Shanghai Jiao Tong University, Shanghai, China, in 2016, and the other in electrical engineering from the University of Notre Dame, South Bend, USA, in 2018. 

He is currently pursuing the Ph.D. degree at Shanghai Jiao Tong University. His research interests are the moment problem and its applications to stochastic filtering, distribution steering, system identification and statistics. 

\begin{wrapfigure}{l}{25mm} 
\includegraphics[width=1in,height=1.25in,clip,keepaspectratio]{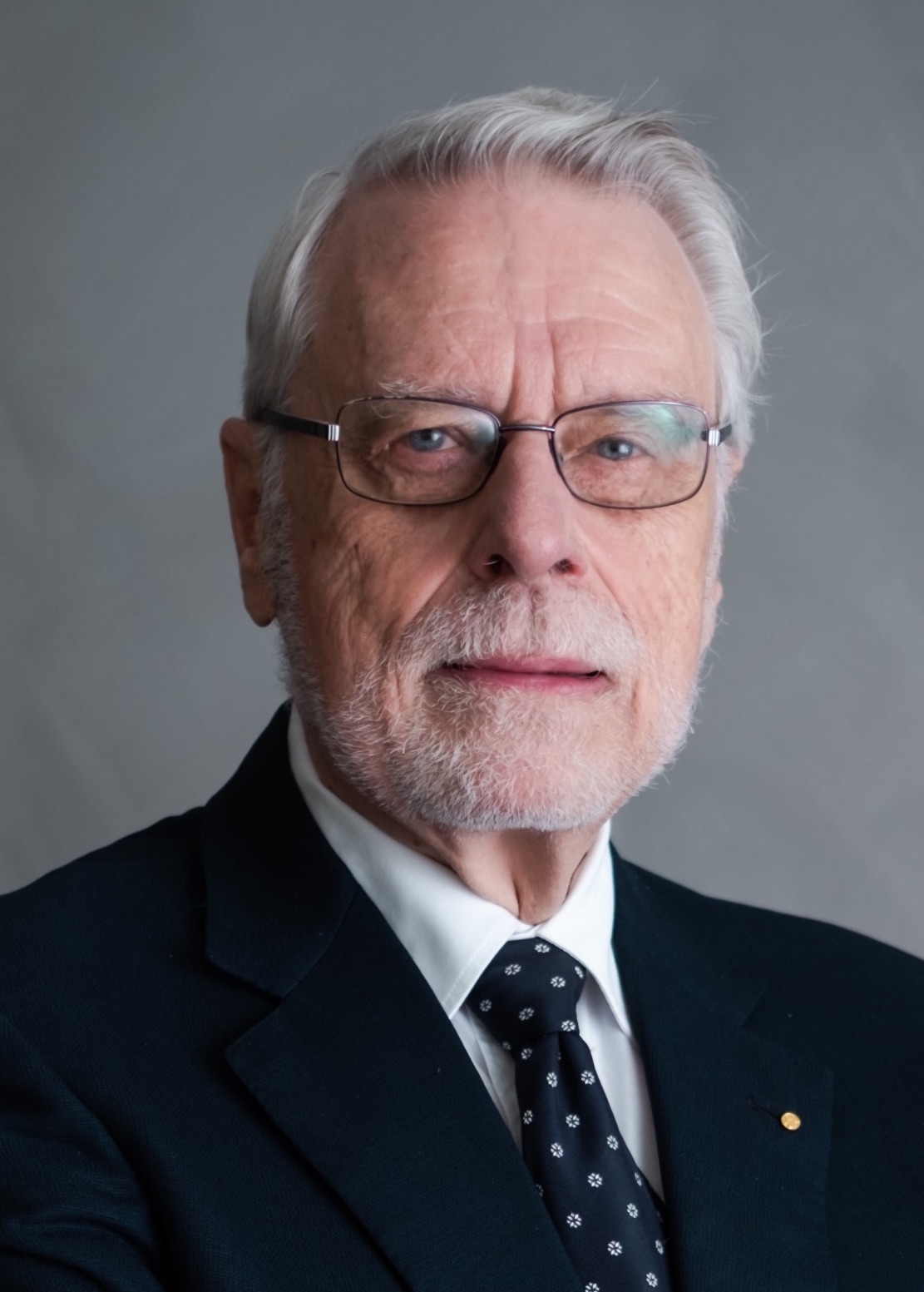}
  \end{wrapfigure}\par
  \textbf{Anders Lindquist} received the Ph.D. degree in Optimization and Systems Theory from the Royal Institute of Technology (KTH), Stockholm, Sweden, in 1972, an honorary doctorate (Doctor Scientiarum Honoris Causa) from Technion (Israel Institute of Technology) in 2010 and Doctor Jubilaris from KTH in 2022.

He is currently a Distinguished Professor at Anhui University, Hefei, China, Professor Emeritus at Shanghai Jiao Tong University, China, and Professor Emeritus at the Royal Institute of Technology (KTH), Stockholm, Sweden. Before that he had a full academic career in the United States, after which he was appointed to the Chair of Optimization and Systems at KTH.

Dr. Lindquist is a Member of the Royal Swedish Academy of Engineering Sciences, a Foreign Member of the Chinese Academy of Sciences, a Foreign Member of the Russian Academy of Natural Sciences (elected 1997), a Member of Academia Europaea (Academy of Europe), an Honorary Member the Hungarian Operations Research Society, a Life Fellow of IEEE, a Fellow of SIAM, and a Fellow of IFAC. He received the 2003 George S. Axelby Outstanding Paper Award, the 2009 Reid Prize in Mathematics from SIAM, and the 2020 IEEE Control Systems Award, the IEEE field award in Systems and Control.

\end{document}